\documentclass[10pt, a4paper]{article}
\usepackage{amsmath, amssymb, amsthm}
\usepackage[top=2cm, bottom=2cm, left=2cm, right=2cm]{geometry}
\usepackage{xcolor}
\usepackage{enumitem}
\usepackage{bbm}
\usepackage[style=numeric, sorting=nty, sortcites=true, backend=bibtex, maxbibnames=10]{biblatex}
\AtEveryBibitem{%
	\clearfield{doi}%
} 
\renewbibmacro*{publisher+location+date}{%
	\ifentrytype{book}
	{\printlist{publisher}%
		\iflistundef{location}
		{\setunit*{\addcomma\space}}
		{\setunit*{\addcomma\space}}%
		\printlist{location}}
	{\printlist{location}%
		\iflistundef{publisher}
		{\setunit*{\addcomma\space}}
		{\setunit*{\addcomma\space}}%
		\printlist{publisher}}
	\setunit*{\addcomma\space}%
	\usebibmacro{date}%
	\newunit} 
\numberwithin{equation}{section}
\newcommand{\pd}[2]{\dfrac{\partial #1}{\partial #2}} 
\newcommand{\td}[2]{\dfrac{\mathrm{d} #1}{\mathrm{d} #2}} 
\newcommand{\pdd}[2]{\frac{\partial #1}{\partial #2}} 
\newcommand{\ii}[4]{\int_{#1}^{#2} #3 \: \mathrm{d}#4} 
\newcommand{\aevv}{\text{ a.e. in } } 
\newcommand{\aevvv}{\text{a.e. in }} 
\newcommand{\aev}{\text{ a.e. } t \geq 0 } 
\newcommand{\dq}[1]{\partial_t^h #1} 
\newcommand{\ene}{\widetilde{\Psi}_\Omega} 
\newcommand{\eneg}{\widetilde{\Psi}_{\Omega, \text{reg}}} 
\newcommand\restr[2]{{
		\left.\kern-\nulldelimiterspace 
		#1 
		\vphantom{\big|} 
		\right|_{#2} 
}} 
\newcommand{\esssup}[2]{\underset{#1}{\text{ess sup}} \: #2}
\DeclareMathOperator{\spn}{span}
\DeclareMathOperator{\dist}{dist}
\DeclareMathOperator{\rg}{rg}

\theoremstyle{remark}
\newtheorem{remark}{Remark}[section]
\theoremstyle{plain}
\newtheorem{proposition}{Proposition}[section]
\newtheorem{lemma}{Lemma}[section]
\newtheorem{theorem}{Theorem}[section]
\theoremstyle{definition}
\newtheorem{definition}{Definition}[section]
\newtheorem*{assumption}{Assumption A}

\begin{document}

\begin{center}
	\large{\textbf{\uppercase{On a system of coupled Cahn-Hilliard equations}}}\\
	\vspace{0.5cm}
	\small{\uppercase{Andrea Di Primio}\footnote[1]{andrea.diprimio$@$polimi.it - Dipartimento di Matematica, Politecnico di Milano, Milano 20133, Italy}},\hspace{0.5cm}
\uppercase{Maurizio Grasselli}\footnote[2]{maurizio.grasselli$@$polimi.it - Dipartimento di Matematica, Politecnico di Milano, Milano 20133, Italy \\ \copyright~2022. Licensed under the Creative Commons CC-BY-NC-ND 4.0 license \url{https://creativecommons.org/licenses/by-nc-nd/4.0/}.}
\end{center}

	
	\begin{abstract}
		\noindent
		We consider a system which consists of a Cahn-Hilliard equation coupled with a Cahn-Hilliard-Oono equation in a bounded domain of $\mathbb{R}^d$, $d = 2, 3$. This system accounts for macrophase and microphase separation in a polymer mixture through two order parameters $u$ and $v$. The free energy of this system is a bivariate interaction potential which contains the mixing entropy of the two order parameters and suitable coupling terms. The equations are endowed with initial conditions and homogeneous Neumann boundary conditions both for $u,v$ and for the corresponding chemical potentials. We first prove that the resulting problem is well posed in a weak sense. Then, in the conserved case, we establish that the weak solution regularizes instantaneously. Furthermore, in two spatial dimensions, we show the strict separation property for $u$ and $v$, namely, they both stay uniformly away from the pure phases $\pm 1$ in finite time. Finally, we investigate the long-time behavior of a finite energy solution showing, in particular, that it converges to a single stationary state.
	\end{abstract}

	\section{Introduction}
	\label{sec:intro}
Cahn-Hilliard type equations are extensively used to model phase separation phenomena which occur in many different contexts (see, for instance, \cite{Miranville19} and references therein). The prototypical example in this regard is represented by phase separation processes taking place in binary alloys, as originally proposed in \cite{CahnHilliard58}. Here we are interested in the theoretical analysis of a system of Cahn-Hilliard equation coupled with a Cahn-Hilliard-Oono equation proposed in \cite{Avalos16} to describe the dynamics of certain polymer blends. In that framework, a mixture consisting of a diblock copolymer and a homopolymer (see \cite{iupac} for detailed definitions) is taken into consideration. From the phenomenological point of view, two distinct, but simultaneous phase separation processes take place. On one hand, the so-called \textit{macrophase} separation occurs between the homopolymer and the copolymer. As a consequence, the diblock copolymer is confined in a region assuming typically a spherical or ellipsoidal shape. On the other hand, the \textit{microphase} separation involves the two blocks of the copolymer, creating regions characterized by the prevalence of one or the other. The two processes generate strikingly regular patterns, which have recently been experimentally investigated as well as analyzed numerically (see, e.g., \cite{Avalos16,Avalos18,Li20,SodiniMartini21} and their references).

In order to introduce the system we indicate by $\Omega \subset \mathbb{R}^d$, $d = 2,3$ a bounded, open, connected and sufficiently smooth domain. Then, we denote with $u(x,t)$ and $v(x,t)$ the relative concentration differences of the phases in the macrophases and the microphase at point $x\in\Omega$ at time $t$, respectively, while the corresponding chemical potentials are instead denoted by $\mu(x,t)$ and $\varphi(x,t)$. Let $T > 0$ be a given final time and $\varepsilon_u^2, \varepsilon_v^2, \sigma$ be three positive real parameters. The system reads as follows (see \cite{Avalos16})
	\begin{equation}
		\label{eq:coupledCH}
		\begin{cases}
			\pd{u}{t} = \Delta \mu & \quad \text{in } \Omega \times (0,T),\\
			\mu = - \varepsilon_u^2 \Delta u + \pd{F}{u}(u,v) & \quad \text{in } \Omega \times (0,T),\\[0.2cm]
			\pd{v}{t} + \sigma\left( v - \dfrac{1}{|\Omega|} \displaystyle\int_\Omega v_0 \: \mathrm{d}x \right) = \Delta \varphi & \quad \text{in } \Omega \times (0,T),\\[0.25cm]
			\varphi = - \varepsilon_v^2 \Delta v + \pd{F}{v}(u,v) & \quad \text{in } \Omega \times (0,T),\\
			\pd{u}{\mathbf{n}} = \pd{v}{\mathbf{n}} = 0 & \quad \text{on } \partial \Omega \times (0,T),\\[0.2cm]
			\pd{\mu}{\mathbf{n}} = \pd{\varphi}{\mathbf{n}} = 0 & \quad \text{on } \partial \Omega \times (0,T),\\[0.2cm]
			u(\cdot, 0) = u_0 & \quad \text{in } \Omega,\\
			v(\cdot, 0) = v_0 & \quad \text{in } \Omega,
		\end{cases}
	\end{equation}
	where $|\Omega|$ denotes the $d$-dimensional Lebesgue measure of $\Omega$ and the bivariate potential density $F$ is the sum of three contributions, namely
	\begin{equation}
		\label{eq:potentialreg}
		F(u,v) = \dfrac{(u^2-1)^2}{4} + \dfrac{(v^2-1)^2}{4} + C(u,v),
	\end{equation}
	where the coupling term $C(u,v)$ is given by
	\begin{equation}
\label{eq:coupling}
C(u,v) = \alpha uv + \beta uv^2 + \gamma u^2v.
\end{equation}
	Here the coupling coefficients $\alpha, \beta, \gamma$ are three given real parameters. In the present work we are interested to provide a theoretical analysis of Problem \eqref{eq:coupledCH} by replacing the double well potentials in \eqref{eq:potentialreg} with the thermodynamically relevant potentials characterized by the mixing entropy densities. More precisely, we consider Problem \eqref{eq:coupledCH} with
	\begin{equation}
		\label{eq:potential}
		F(u,v) = S(u;\theta_u,\theta_{0,u}) + S(v;\theta_v,\theta_{0,v}) + C(u,v),
	\end{equation}
	where
	\begin{equation}
    \label{eq:mixentropy}
    S(r, \theta_r, \theta_{0,r}) = \dfrac{\theta_r}{2}[(1+r)\log(1+r) + (1-r)\log(1-r)] - \dfrac{\theta_{0,r}}{2}r^2 , \quad r \in (-1,1),
    \end{equation}
	with $0 < \theta_r < \theta_{0,r}$ and $r = u$ or $r=v$. Here $\theta_r$ and $\theta_{0,r}$ represent the absolute temperature and the critical temperature under which the separation processes take place, respectively. We recall that $S$ is known as singular potential (or Flory--Huggins potential, see \cite{Flory42,Huggins41}, cf. also \cite{Glotzer}). We point out that the regular double well potentials in \eqref{eq:potentialreg} are just
convenient approximations of $S$ but they do not ensure that $u$ and $v$ take their values in the physical range $[-1,1]$.

It is worth recalling that $\varepsilon_u^2, \varepsilon_v^2, \sigma$, as well as the coupling coefficients, have a physical interpretation in the framework of polymer blends. For instance, the quantities $\varepsilon_u$ and $\varepsilon_v$ are proportional to the thickness of the propagating fronts of each component, and are therefore linked to the rapidity of variation of $u$ and $v$ in the interface region (see \cite{Avalos16} for the details).


		
Problem \eqref{eq:coupledCH} entails the conservation of the total mass of both the order parameters. Indeed, setting
		\[ \overline{f}:= \dfrac{1}{|\Omega|} \int_\Omega f \: \mathrm{d}x, \]
		for any Lebesgue-integrable function $f$, then we get
		\begin{equation}
			\label{eq:avg}
			\begin{cases}
				\overline{u}(t) = \overline{u_0},\\ \overline{v}(t) = \overline{v_0},\\
			\end{cases}
		\end{equation}
		for any $t \geq 0$. However, this might not always be the case. Indeed, if we consider the Cahn-Hilliard-Oono equation for $v$ in the following general form (see, e.g., \cite{GGM} and references therein)
		\[ \pd{v}{t} + \sigma\left( v - c \right) = \Delta \varphi \qquad \text{in } \Omega \times (0,T), \]
		for some prescribed $c \in (-1,1)$, one obtains
		\[ \overline{v}(t) = c + e^{-\sigma t}(\overline{v_0} - c), \]
		implying that two possible scenarios may arise. In the conserved case, the quantity $\overline{v}$ is constant and equal to $c = \overline{v_0}$, whereas in the so-called \textit{off-critical} case, i.e. $c \neq \overline{v_0}$, $\overline{v}(t)$ converges exponentially fast to $c$ as $t$ approaches infinity. As we shall see, this is not a small detail from the theoretical viewpoint (cf. \cite{GGM}).

We also note that the conserved case can be seen as the gradient flow of the free energy
		\begin{equation}
			\label{eq:ohtakawasaki}
			\Psi_\Omega(u, \nabla u, v, \nabla v) = \int_\Omega \left(\varepsilon_u^2 \dfrac{|\nabla u|^2}{2} + \varepsilon_v^2 \dfrac{|\nabla v|^2}{2} + F(u,v)\right) \: \mathrm{d}x  + \sigma\int_\Omega\int_\Omega (v(x)-\overline{v_0})G(x,y)(v(y)-\overline{v_0}) \: \mathrm{d}y\mathrm{d}x,
		\end{equation}
		which is known as Ohta-Kawasaki functional (see, for instance, \cite{Choksi03, dcp1, dcp2} and references therein for regular potentials). In \eqref{eq:ohtakawasaki}, $G$ denotes the Green function associated to the negative Laplace operator with homogeneous Neumann boundary conditions.
	\noindent
We recall that the Cahn-Hilliard-Oono equation with singular potential has recently been analyzed in \cite{GGM} while its coupling with the Navier-Stokes system has been studied in \cite{Miranville2} (see also their references for the regular potential case). Instead, only numerical simulations are available so far for Problem \eqref{eq:coupledCH}. Our goal is
to extend the analysis done in \cite{GGM} to the present problem. As we shall see, this is not a straightforward task because of the coupling term \eqref{eq:coupling}.
\\[\baselineskip]
	\textbf{Plan of the paper.} In Section \ref{sec:notation}, we introduce some notation and the functional setting. Section \ref{sec:wellposed} is devoted to introduce a weak formulation of our problem in the off-critical case and to state its well-posedness  whose proof is given in Section \ref{sec:wpproof}. The regularization properties in the conserved case are analyzed in Section \ref{sec:regularity}, while Section \ref{sec:separation} is devoted to establish, in the conserved case, the strict separation property of both the macrophase and the microphase in dimension two. In Section \ref{sec:equilibrium} we analyze the longtime behavior of weak solutions in the conserved case. In particular, we show that any weak solution converges to a single stationary state.
	
	\section{Notation and functional setting}
	\label{sec:notation}
	Throughout all this work, given any pair of positive integers $k,p$, we denote by $W^{k,p}(\Omega)$ the Sobolev space of $L^p(\Omega)$ functions with distributional derivatives of order less or equal to $k$ also in $L^p(\Omega)$. This space is endowed with the classical norm $\| \cdot \|_{W^{k,p}(\Omega)}$. For any choice of $k \in \mathbb{N}$, $H^k(\Omega):=W^{k,2}(\Omega)$ is a Hilbert space with respect to the scalar product
	\[ (f,g)_{H^k} := \sum_{|s| \leq k} (D^sf, D^sg),\]
	for any $f,g \in H^k(\Omega)$. We recall the Hilbert triplet
	\[
	V := H^1(\Omega) \hookrightarrow H := L^2(\Omega) \hookrightarrow V^* := H^1(\Omega)^*,
	\]
	with dense, continuous and compact injections (in both two and three spatial dimensions). Here $H^1(\Omega)^*$ denotes the topological dual space of $H^1(\Omega)$. The three spaces $V, H, V^*$ are endowed with the norms $\| \cdot \|_V, \| \cdot \|, \| \cdot \|_{V^*}$, respectively. In particular, $\| \cdot \|$ denotes the classical $L^2$-norm (possibly for functions taking values in $\mathbb{R}^d$), whereas
	\[ \|u\|^2_V = \|u\|^2 + \| \nabla u \|^2 \]
	for every $u \in V$. Finally, $\| \cdot\|_{V^*}$ is the standard operator norm in a dual space. From here onwards, the scalar products inducing said norms are denoted as $(\cdot, \cdot)_\cdot$, accordingly. The duality between a (real) Banach space $X$ and its topological dual $X^*$ is denoted by $\langle \cdot, \cdot \rangle$. We now recall some well-known and useful results. Let us introduce the spaces
	\[
	V_0 := \{ u \in V: \overline{u} = 0\}, \qquad \quad V_0^* := \left\{ L \in V^*: \overline{L} := \dfrac{1}{|\Omega|}\langle L, 1 \rangle = 0\right\},
	\]
	where $\langle \cdot, \cdot \rangle$ denotes the pairing between $V^*$ and $V$. Let us consider the linear operator
	\begin{equation*}
		A:  V \to V^*,\qquad u \mapsto \left( v \mapsto \ii{\Omega}{}{\nabla u \cdot \nabla v}{x} \right),
	\end{equation*}
	whose restriction to $V_0$ is an isomorphism between $V_0$ and its topological dual $V_0^*$.
\noindent
The inverse operator $\mathcal{N} := A^{-1}$ is well defined. This operator, by definition, satisfies
\begin{equation}
	\label{eq:inverseAN}
	\mathcal{N}Au = u, \qquad A\mathcal{N}L = L \qquad \forall u \in V_0,\; \forall  L \in V_0^*,
\end{equation}
so that $\mathcal{N}$ is the inverse of the negative Laplace operator with homogeneous Neumann conditions. The following result is useful and straightforward to prove.
\begin{proposition}
	\label{prop:propertiesAN}
	Let $A$ and $\mathcal{N}$ be defined as above. Then, the following identities hold
	\begin{enumerate}[label=(\roman*)]
	\item $ \langle Au, \mathcal{N}L \rangle = \langle L,u \rangle, \quad \forall \: u \in V_0, \: L \in V^*_0$;
	\item $ \langle L_1, \mathcal{N}L_2 \rangle = (\nabla (\mathcal{N}L_1), \nabla (\mathcal{N}L_2)), \quad \forall \: L_1, L_2 \in V_0^*.$
	\end{enumerate}	
\end{proposition}
\noindent
Concerning the choice of suitable norms on $V_0$ and $V_0^*$, we define
\[ \| L \|_* := \|\nabla(\mathcal{N}L) \| = \sqrt{\langle L, \mathcal{N}L \rangle} \]
on $V_0^*$, where the second equality is due to Proposition \ref{prop:propertiesAN}-(ii), and
\[ \| L \|_{-1}^2 := \| L - \overline{L} \|^2_* + |\overline{L}|^2 \]
on $V^*$. Finally, we state the following
\begin{proposition}
	\label{prop:normeq}
	The norm $\| \cdot \|_*$ is an equivalent norm in $V_0^*$, while the norm $\| \cdot \|_{-1}$ is an equivalent norm in $V^*$ (with respect to the usual operator norms).
\end{proposition}

\section{Weak formulation and well-posedness}
\label{sec:wellposed}
In this section we first introduce the weak formulation of \eqref{eq:coupledCH} in the more general off-critical case, namely, the third equation is replaced by
\[
\pd{v}{t} + \sigma\left( v - c \right) = \Delta \varphi \qquad \text{in } \Omega \times (0,T),
\]
for some known constant $c \in (-1,1)$. More precisely, we consider the following (formal) initial and boundary value problem
\begin{equation}
	\label{eq:coupledCHoc}
	\begin{cases}
		\pd{u}{t} = \Delta \mu & \quad \text{in } \Omega \times (0,T),\\
		\mu = - \varepsilon_u^2 \Delta u + \pd{F}{u}(u,v) & \quad \text{in } \Omega \times (0,T)\\[0.2cm]
		\pd{v}{t} + \sigma\left( v - c \right) = \Delta \varphi & \quad \text{in } \Omega \times (0,T),\\[0.25cm]
		\varphi = - \varepsilon_v^2 \Delta v + \pd{F}{v}(u,v) & \quad \text{in } \Omega \times (0,T),\\
		\pd{u}{\mathbf{n}} = \pd{v}{\mathbf{n}} = 0 & \quad \text{on } \partial \Omega \times (0,T),\\[0.2cm]
		\pd{\mu}{\mathbf{n}} = \pd{\varphi}{\mathbf{n}} = 0 & \quad \text{on } \partial \Omega \times (0,T),\\[0.2cm]
		u(\cdot, 0) = u_0 & \quad \text{in } \Omega,\\
		v(\cdot, 0) = v_0 & \quad \text{in } \Omega.
	\end{cases}
\end{equation}
The main goal of this section is to state and prove the well-posedness of the weak formulation of \eqref{eq:coupledCHoc} which is given by
\begin{definition}
	\label{def:solution}
	Let $u_0, v_0 \in V$ be such that $F(u_0, v_0) \in L^1(\Omega)$ and $\overline{u_0}, \overline{v_0} \in (-1,1)$. Let $c \in (-1,1)$. A weak solution to Problem \eqref{eq:coupledCHoc} is a pair $(u,v)$ enjoying the following properties:
	\begin{enumerate}[label=(\roman*)]
		\item $u \in L^{\infty}([0,T]; V) \cap L^2([0,T]; H^2(\Omega)) \cap L^{\infty}(\Omega \times (0,T))$;
		\item $v \in L^{\infty}([0,T]; V) \cap L^2([0,T]; H^2(\Omega)) \cap L^{\infty}(\Omega \times (0,T))$;
		\item $\frac{\mathrm{d}u}{\mathrm{d}t} \in L^2([0,T]; V^*)$;
		\item $\frac{\mathrm{d}v}{\mathrm{d}t} \in L^2([0,T]; V^*)$;
		\item $\mu = -\varepsilon_u^2 \Delta u + \frac{\mathrm{d}F}{\mathrm{d}u} \in L^2([0,T]; V)$;
		\item $\varphi = -\varepsilon_v^2 \Delta v + \frac{\mathrm{d}F}{\mathrm{d}v} \in L^2([0,T]; V)$;
		\item $|u(\mathbf{x}, t)| < 1$ for a.a. $(\mathbf{x},t) \in \Omega \times (0,T)$;
		\item $|v(\mathbf{x}, t)| < 1$ for a.a. $(\mathbf{x},t) \in \Omega \times (0,T)$;
		\item $(u,v)$ solves the system
		\[
     	\begin{cases}
		\label{eq:weak2}	
			\left\langle \pd{u}{t}, s \right\rangle + (\nabla \mu, \nabla s) = 0  & \qquad \forall \: s \in V, \: \aevv (0,T),\\
			\left\langle \pd{v}{t}, w \right\rangle + \sigma(v-c, w) + (\nabla \varphi, \nabla w) = 0 & \qquad \forall \: w \in V, \: \aevv (0,T);
		\end{cases}
		\]
		\item $\pdd{u}{\mathbf{n}} = 0$ a.e. in $\partial\Omega \times (0,T)$;
		\item $\pdd{v}{\mathbf{n}} = 0$ a.e. in $\partial\Omega \times (0,T)$;
        \item $u(0)=u_0$ a.e. in $\Omega$;
        \item $v(0)=v_0$ a.e. in $\Omega$.
	\end{enumerate}
\end{definition}
\begin{remark}
	\label{rem:reg1}
	From Definition \ref{def:solution}, thanks to the fact that $u,v \in L^2([0,T];V)$ and their time derivatives belong to $L^2([0,T];V^*)$, we directly infer that $u,v \in \mathcal{C}^0([0,T];H)$.
\end{remark}
\begin{remark}
	\label{rem:reg3}
	Let $T > 0$ be arbitrary. On account of the $L^\infty(\Omega \times (0,T))$-regularity of each phase, it holds that $u,v \in L^\infty([0,T]; L^p(\Omega))$ for each $p \geq 1$. In particular, the function mapping $t \mapsto \| u(t) \|_{L^\infty}$ (same for $v$) is measurable and essentially bounded (see \cite[Rem. 3.3]{GGG}).
\end{remark}
\noindent
Let us define now the energy functional
\begin{equation}
\label{eq:ENERID}
\Psi_\Omega(u,v)  = \dfrac{\varepsilon_u^2}{2} \| \nabla u \|^2 + \dfrac{\varepsilon_v^2}{2} \| \nabla v \|^2 + \ii{\Omega}{}{F(u,v)}{x}.
\end{equation}
The well-posedness result is
\begin{theorem}
	\label{th:wellposed}
	Let $u_0,v_0 \in V$ be such that $F(u_0,v_0) \in L^1(\Omega)$ and $\overline{u_0}, \overline{v_0} \in (-1,1)$. Let $c \in (-1,1)$ be given.
Then, Problem \eqref{eq:coupledCHoc} has a unique finite energy solution. Moreover, the following energy inequality holds
	\begin{equation}
		 \label{eq:dissipative}
		 \Psi_\Omega(u(t),v(t)) + \dfrac{1}{2}\int_t^{t+1} \left(\|\nabla \mu(\tau) \|^2 + \| \nabla \varphi(\tau) \|^2 \right)\: \mathrm{d}\tau \leq
         \Psi_\Omega(u_{0},v_{0})e^{-\sigma t} + C
	\end{equation}
for any $t \geq 0$, where $C$ is a positive constant depending on all the parameters of the problem. Furthermore, given $R \geq 0$, $T > 0$ and $m \in (0,1]$ such that $|c| \leq m$, there exists a constant $K = K(m,R,T)$ such that, for any solutions $(u_1, v_1), (u_2,v_2)$ on $[0,T]$ originating from the initial conditions $(u_{01}, v_{01}), (u_{02}, v_{02})$ satisfying $\Psi_\Omega(u_{0i}, v_{0i}) \leq R$ and $|\overline{u_{0i}}|, |\overline{v_{0i}}| \leq m$ for $i = 1,2$, the continuous dependence estimate
	\begin{multline}
		\| u_1(t) - u_2(t) \|_{V^*} + \| v_1(t) - v_2(t) \|_{V^*} + \left( \int_0^T \|u_1(t)-u_2(t)\|^2_{V} \: \mathrm{d}t \right)^\frac{1}{2}+ \left( \int_0^T \|v_1(t)-v_2(t)\|^2_{V} \: \mathrm{d}t \right)^\frac{1}{2}\\ \leq K(\|u_{01} - u_{02}\|_{V^*} + \|v_{01} - v_{02}\|_{V^*} + |\overline{u_{01}}-\overline{u_{02}}|^{\frac{1}{2}} + |\overline{v_{01}}-\overline{v_{02}}|^{\frac{1}{2}})
\label{contdep}
	\end{multline}
	holds for every $t \in [0,T]$ and entails the uniqueness of a weak solution.
\end{theorem}

\begin{remark}
Elliptic arguments yield higher regularity for $u$ and $v$. Indeed, on account of \cite[(3.11)]{GGW}, we can also deduce $u, v \in L^4([0,T]; H^2(\Omega))$. Also, recalling \cite[Thm. 6]{Abels}, we can prove that $u, v \in L^2([0,T]; W^{2,r}(\Omega))$, where $r=6$ if $d=3$ or $r\in (1,\infty)$ if $d=2$.
\end{remark}

\section{Proof of Theorem \ref{th:wellposed}}
\label{sec:wpproof}
The proof is split into four steps. First we introduce a convenient approximation of the potential $F$ which is crucial in order to establish the existence
of a weak solution. This is obtained in the second step through a suitable Galerkin scheme. The final two steps are devoted to prove the global energy inequality  \eqref{eq:dissipative} and the continuous dependence estimate \eqref{contdep}.

\subsection{Approximating the bivariate potential}
\label{subsec:apprF}
Recalling \eqref{eq:mixentropy}, we set
\begin{equation} \label{eq:singularpart}
	\hat{S}(r;\theta_r) = \dfrac{\theta_r}{2}\left[(1+r)\log(1+r) + (1-r)\log(1-r)\right], \qquad r \in (-1,1).
\end{equation}
We point out that $\hat{S}$ is meant to be extended by (right or left) continuity at the endpoints and then over the whole real line with value $\hat{S}(r) = +\infty$ whenever $|r|> 1$. It is well known that the function $\hat{S}$ has the following elementary properties:
\begin{enumerate}[label = (\roman*)]
	\label{page:test}
	\item $\hat{S}$ is real analytic in $(-1,1)$, and in particular belongs to $\mathcal{C}^4((-1,1))$, where $\mathcal{C}^p(I)$ denotes the set of (classically) $k$-times-continuously differentiable functions over an interval $I$ when $p > 0$, and the set of continuous functions over $I$ when $p = 0$;
    	\item there exists a constant $k > 0$ such that $\hat{S}^{(4)}(u; \theta_u)$ and $\hat{S}^{(4)}(v; \theta_v)$  are non-decreasing in $(-1, -1+k)$ and non-increasing in $(1-k, 1)$;
	\item there holds
	\[ \lim_{r \to -1^+} \hat{S}'(r) = -\infty; \qquad \lim_{r \to 1^-} \hat{S}'(r) = +\infty; \]
	\item there holds
\[\hat{S}''(r; \theta_u) \geq \theta_u > 0; \qquad \hat{S}''(r; \theta_v) \geq \theta_v > 0, \qquad \forall \: r \in(-1,1); \]
	\item there exists $c > 0$ and $\epsilon_0 > 0$ such that
	\[
	\hat{S}^{(4)}(r; \theta_u) \geq c; \qquad \hat{S}^{(4)}(r; \theta_v)  \geq c, \qquad \forall \: r \in (-1, -1 + \epsilon_0] \cup [1 - \epsilon_0, 1);
	\]
	\item there exists $\epsilon_1 > 0$ such that, for each $k = 0, 1, 2, 3, 4,$ and for each $j = 0,1$,
	\[ \begin{cases}
		\hat{S}^{(k)}(r; \theta_r) \geq 0 & \quad \forall \: r \in [1 - \epsilon_1, 1),\\
		\hat{S}^{(2j+2)}(r; \theta_r) \geq 0 & \quad \forall \: r \in (-1, -1+\epsilon_1],\\
		\hat{S}^{(2j+1)}(r; \theta_r) \leq 0 & \quad \forall \: r \in (-1, -1+\epsilon_1],\\
	\end{cases}\]
	in both cases $r = u$ and $r = v$, where $\hat{S}^{(0)} := \hat{S}$.
\end{enumerate}
\color{black}
\begin{remark}
\label{genpot}
In the sequel, $\hat{S}(r;\theta_r)$ will denote any function satisfying the above mentioned properties. Of course, \eqref{eq:singularpart} is an admissible choice.
\end{remark}
\color{black}
The aim of this first step is to introduce a suitable \textit{regular} (i.e. with no singularities over the whole plane $\mathbb{R}^2$) approximation of $F$ (i.e. of $\hat{S}$), dependent on a positive real (small) parameter $\delta$ in such a way that the original potential is recovered in the limit $\delta \to 0^+$. To this end, we introduce a family of regular potentials based upon the fourth-order Taylor expansion of $\hat{S}$ (see \cite{FG12}).
Fixed any sufficiently small $\delta \in (0,1)$, let $\hat{S}_\delta: \mathbb{R} \to \mathbb{R}$ be a globally defined approximation of the singular part of the function $S$ given by
\begin{equation}
	\label{eq:taylor}
	\hat{S}_\delta(r) =
	\begin{cases}
		\displaystyle\sum_{i=0}^{4} \dfrac{\hat{S}^{(i)}(-1+\delta)}{i!}[r-(-1+\delta)]^i & \quad \text{if } r \leq -1 + \delta,\\[0.2cm]
		\hat{S}(r) & \quad \text{if } |r| \leq 1-\delta,\\
		\displaystyle\sum_{i=0}^{4} \dfrac{\hat{S}^{(i)}(1-\delta)}{i!}[r-(1-\delta)]^i& \quad \text{if } r \geq 1 - \delta.
	\end{cases}
\end{equation}
Accordingly, we set
\begin{equation}
		S_\delta(u;\theta_u,\theta_{0,u}) := \hat{S}_\delta(u;\theta_u) - \dfrac{\theta_{0,u}}{2}u^2,\qquad
		S_\delta(v;\theta_v,\theta_{0,v}) := \hat{S}_\delta(v;\theta_v) - \dfrac{\theta_{0,v}}{2}v^2,\\
\end{equation}
and
\begin{equation}
	\label{eq:potdelta}
	F_\delta(u,v) := S_\delta(u;\theta_u,\theta_{0,u}) + S_\delta(v;\theta_v,\theta_{0,v}) + C(u,v).
\end{equation}
From now on the dependence on the absolute and the critical temperatures in $S$, $\hat{S}$ and their regular approximations will be omitted.
Here below we state and prove a result on the coercivity of $F_\delta$ which will be helpful in the next subsection.
\begin{proposition}
	\label{prop:fdelta}
	$F_\delta \in \mathcal{C}^4(\mathbb{R}^2)$ for any sufficiently small $\delta\in (0,1)$. Furthermore, there exists $\delta_0 \in (0,1)$ and two positive constants $k_1, k_2$ independent of $\delta$ such that
	\[ F_\delta(u,v) \geq k_1(u^4 + v^4) - k_2 \qquad \forall \: u,v \in \mathbb{R}, \: \delta \in (0, \delta_0]. \]
\end{proposition}
\begin{proof}
We slightly adapt the proof of \cite[Lemma 1]{FG12}. Without loss of generality, let us first consider $u \geq 0$. With reference to properties (v) and (vi), fix $\delta < \delta_0 := \min(\epsilon_0, \epsilon_1)$. Then, given a value of $u$, one and only one of the following cases applies.
	\begin{itemize}[label=$*$]
		\item \fbox{$0 \leq u \leq 1-\delta$} In this range, $\hat{S}_\delta(u) = \hat{S}(u) \geq 0$. Given any $k_1, k_2 > 0$, then one has
		\[ k_1 u^4 - k_2 \leq k_1 - k_2 \]
		thus, provided that $k_2 \geq k_1$, the right hand side is negative, yielding
		\begin{equation}
			\label{eq:bound1}
			\hat{S}_\delta(u) \geq k_1 u^4 - k_2.
		\end{equation}
		\item \fbox{$1-\delta < u < 1$} Owing to properties (v) and (vi), one has
		\[
		\hat{S}_\delta(u) \geq \dfrac{1}{24}\hat{S}^{(4)}(1-\delta)[u - (1-\delta)]^{4} \geq \dfrac{c}{24}[u - (1-\delta)]^{4} \geq 0,
		\]
		and thus one argues similarly to get \eqref{eq:bound1}.\\
		\item \fbox{$u \geq 1$} From properties (v) and (vi) one easily proves
		\[\hat{S}_\delta(u) \geq \dfrac{1}{24}\hat{S}^{(4)}(1-\delta)[u - (1-\delta)]^{4} \geq \dfrac{c}{24}(u-1)^{4} \geq \dfrac{c}{48}u^{4} - k,\]
		so that we get $k_1 = \frac{c}{48}$ and $k_2 \geq k = k(c)$.
		
	\end{itemize}
	As far as the case $u < 0$ is concerned, a very similar reasoning can be carried out, so that, in conclusion,
	\begin{equation}
		\label{eq:prop1}
		\hat{S}_\delta(u) \geq k_1u^4 - k_2 \qquad \forall \: u \in \mathbb{R},
	\end{equation}
	for a fixed $k_1$ and $k_2 \geq k = k(c)$, provided that $\delta$ is sufficiently small. This, of course, lets us deduce that
	\begin{equation}
		\label{eq:prop2}
		\hat{S}_\delta(u) + \hat{S}_\delta(v) \geq k_1(u^4+v^4) - k_2 \qquad \forall \: u,v \in \mathbb{R},
	\end{equation}
	for any sufficiently small $\delta$.
	Finally, we consider the polynomial term. Owing to the elementary inequality $x \geq -|x|$ and the Young inequality,
	\begin{multline*}
		C(u,v) - \dfrac{\theta_{0,u}}{2}u^2 - \dfrac{\theta_{0,v}}{2}v^2 = \alpha uv + \beta uv^2 + \gamma u^2v - \dfrac{\theta_{0,u}}{2}u^2 - \dfrac{\theta_{0,v}}{2}v^2 \geq\\ \geq - \dfrac{\theta_{0,u}}{2}u^2 - \dfrac{\theta_{0,v}}{2}v^2 -|\alpha|\left(\dfrac{u^2}{2} + \dfrac{v^2}{2}\right) -\left(\dfrac{|\beta|}{3}+\dfrac{2|\gamma|}{3}\right)|u|^3 -\left(\dfrac{2|\beta|}{3}+\dfrac{|\gamma|}{3}\right) |v|^3.
	\end{multline*}
	Let $\lambda > 0$ be arbitrary. Then, there exists a constant $K = K(\lambda)$ such that
	\begin{equation}
		\label{eq:prop3}
		C(u,v) - \dfrac{\theta_{0,u}}{2}u^2 - \dfrac{\theta_{0,v}}{2}v^2 \geq - \lambda(u^4 + v^4) - K(\lambda) \qquad \forall \: u,v \in \mathbb{R},
	\end{equation}
	where $K(\lambda)$ is a constant depending only on the choice of $\lambda$. In particular, if $\lambda < k_1$, then, up to redefinition of constants, then \eqref{eq:prop2} and \eqref{eq:prop3} imply the thesis.	
\end{proof}
\begin{remark}
Proposition \ref{prop:fdelta} shows why \eqref{eq:taylor} is particularly convenient in this context. Indeed if we chose a more elegant approximation like, for instance, the one in \cite{GGG} then we were not able to control the growth of the interaction polynomial $C(u,v)$.
\end{remark}

\subsection{Existence}
\label{subsec:exist}
Let $\delta \in (0,1)$ be sufficiently small. Consider an approximation of Problem \eqref{eq:coupledCHoc} with $F_\delta$ defined by \eqref{eq:potdelta} in place of $F$, namely
	\begin{equation}
		\label{eq:step2}
		\begin{cases}
			\pd{u}{t} = \Delta \mu & \quad \text{in } \Omega \times (0,T),\\
			\mu = -\varepsilon_u^2\Delta u + \dfrac{\partial  F_\delta}{\partial u}(u,v) & \quad \text{in } \Omega \times (0,T),\\
			\pd{v}{t} + \sigma (v - c)  = \Delta \varphi & \quad \text{in } \Omega \times (0,T),\\
			\varphi = -\varepsilon_v^2 \Delta v + \dfrac{\partial F_\delta}{\partial v}(u,v) & \quad \text{in } \Omega \times (0,T),\\
			\pd{u}{\mathbf{n}} = \pd{v}{\mathbf{n}} = 0 & \quad \text{on } \partial\Omega \times (0,T),\\[0.2cm]
			\pd{\mu}{\mathbf{n}} = \pd{\varphi}{\mathbf{n}} = 0 & \quad \text{on } \partial\Omega \times (0,T),\\
			u(\cdot, 0) = u_0 & \quad \text{in } \Omega,\\
			v(\cdot, 0) = v_0 & \quad \text{in } \Omega.
		\end{cases}
	\end{equation}
	The weak formulation of Problem \eqref{eq:step2} is similar to the one for Problem \eqref{eq:coupledCHoc}, provided that $\mu$ and $\varphi$ are now computed using $F_\delta$ instead of $F$. In particular, we have
	\[
	\begin{cases}
		\left\langle \pd{u}{t}, s \right\rangle + (\nabla \mu, \nabla s) = 0  & \quad \forall \: s \in V, \aevv (0,T),\\
		\left\langle \pd{v}{t}, w \right\rangle + \sigma(v-c, s) + (\nabla \varphi, \nabla w) = 0 & \quad  \forall \: w \in V, \aevv (0,T).
	\end{cases}
	\]
We now establish  the existence of a weak solution to this approximating problem by means of a Galerkin scheme.

In the following, for any $m \in \mathbb{N}$, $m \geq 1$, we denote with $(\eta_m, w_m) \in \mathbb{R} \times H$ the (countably many) eigencouples of the Neumann--Laplace operator, namely the relation
	\[ -\Delta w_m = \eta_m w_m\]
	holds for any positive integer $m$. We recall that the set of eigenvectors is an orthonormal basis in $H$ and an orthogonal basis in $V$. Let $m \in \mathbb{N}$, $m \geq 1$ and set $W_m := \spn(w_1, ..., w_m)$. Consider the projection of the weak form of Problem \eqref{eq:step2} on $W_m$, namely find $u_m, v_m:[0,T] \to W_m$ so that
	\begin{equation}
		\label{eq:step3}
		\begin{cases}
			\left\langle \pd{u_m}{t}, s \right\rangle + (\nabla \mu_m, \nabla s) = 0  & \quad \forall \: s \in W_m, \aevv (0,T),\\
			\mu_m = \Pi_m\left(-\varepsilon_u^2 \Delta u_m + \pd{F_\delta}{u}(u_m, v_m)\right) & \quad \aevvv (0,T),\\
			\left\langle \pd{v_m}{t}, w \right\rangle + \sigma(v_m-c, w) + (\nabla \varphi_m, \nabla w) = 0 & \quad \forall \: w \in W_m, \aevv (0,T),\\
			\varphi_m = \Pi_m\left(-\varepsilon_v^2 \Delta v_m + \pd{F_\delta}{v}(u_m, v_m)\right) & \quad \aevvv (0,T),\\
			u_m(0) = u_{0,m} & \quad \text{in } \Omega,\\
			v_m(0) = v_{0,m} & \quad \text{in } \Omega,
		\end{cases}
	\end{equation}
	where $\Pi_m: H \to W_m$ is the projector onto the finite-dimensional space $W_m$. We denote the coordinates of $u_m$ and $v_m$ with respect to the chosen basis of $W_m$ as
	\[ u_m(t) = \sum_{k = 1}^m u_{m,k}(t)w_k; \qquad v_m(t) = \sum_{k = 1}^m v_{m,k}(t)w_k,\]
	where $u_{m,k}$ and $v_{m,k}$ are real-valued functions from $[0,T]$ for any value of $m$ and $1\leq k \leq m$. As usual in Galerkin schemes, we now take in \eqref{eq:step3} $s = w_i$ and $w = w_j$, as $i,j$ vary from $1$ to $m$, yielding, after some manipulations,
	\begin{equation*}
		\begin{cases}
			\td{u_{m,i}}{t} + \varepsilon_u^2 \eta_i^2 u_{m,i} + \eta_i \left( \Pi_m \left(\pd{F_\delta}{u}(u_m,v_m) \right), w_i\right) = 0 & \quad \forall \: i = 1,...,m, \aevv (0,T),\\[0.3cm]
			\td{v_{m,j}}{t} + (\sigma + \varepsilon_v^2 \eta_j^2) v_{m,j} - \sigma c \: \overline{w_j} + \eta_j \left(\Pi_m \left(\pd{F_\delta}{v}(u_m,v_m) \right), w_j\right) = 0 & \quad \forall \: j = 1,...,m, \aevv (0,T),\\
			u_m(0) = u_{0,m} & \quad \text{in } \Omega,\\
			v_m(0) = v_{0,m} & \quad \text{in } \Omega,
		\end{cases}
	\end{equation*}
	namely a Cauchy initial value problem consisting of $2m$ ordinary differential equations in the unknowns $u_{m,k}, v_{m,k}$, for $k = 1,...,m$. In order to have a better understanding of the differential problem above, we restate it in vectorial form. Let $\mathbf{U}_m$ and $\mathbf{V}_m$ denote the vectors of functions whose components are $u_{m,k}$ and $v_{m,k}$, respectively. Let also $\mathbf{U}_0, \mathbf{V}_0$ denote the real coordinates of $u_{0,m}$ and $v_{0,m}$ with respect to the eigenvector basis, accordingly. Then, we have
	\begin{equation}
		\label{eq:ode}
		\begin{cases}
			\td{}{t}
			\begin{bmatrix}
				\mathbf{U}_m\\
				\mathbf{V}_m
			\end{bmatrix}
			= \tilde{\mathbf{F}}_\delta\left(\begin{bmatrix}
				\mathbf{U}_m\\
				\mathbf{V}_m
			\end{bmatrix} \right) & \quad \aevv (0,T),\\[0.35cm]
			\begin{bmatrix}
				\mathbf{U}_m\\
				\mathbf{V}_m
			\end{bmatrix}(0) = \begin{bmatrix}
				\mathbf{U}_0\\
				\mathbf{V}_0
			\end{bmatrix},
		\end{cases}
	\end{equation}
	where the function $\tilde{\mathbf{F}}_\delta : \mathbb{R}^m \times \mathbb{R}^m \cong \mathbb{R}^{2m} \to \mathbb{R}^{2m}$ has the form
	\[ \tilde{\mathbf{F}}_\delta\left(\begin{bmatrix}
		\mathbf{x}\\
		\mathbf{y}
	\end{bmatrix} \right) = \mathbf{L}\begin{bmatrix}
		\mathbf{x}\\
		\mathbf{y}
	\end{bmatrix} + \mathbf{N}_\delta\left(\begin{bmatrix}
		\mathbf{x}\\
		\mathbf{y}
	\end{bmatrix}\right) + \mathbf{b}\]
	and in turn, $\textbf{L}$ is a diagonal real $2m \times 2m$ matrix, while $\textbf{N}_\delta: \mathbb{R}^{2m} \to \mathbb{R}^{2m}$ is a nonlinear function depending on the potential approximation parameter $\delta$ and $\mathbf{b} \in \mathbb{R}^{2m}$ is a constant real vector. In light of the Cauchy-Lipschitz theorem, system \eqref{eq:ode} has a unique solution
$(\mathbf{U}_m,\mathbf{V}_m) \in C^1([0,T^*_m]; \mathbb{R}^{2m})$ on the maximal interval $[0,T^*_m]$ where $T^*_m\in (0,T]$. We now show that the approximating (local) solution satisfies an energy inequality itself. Let us consider the approximated energy functional
	\[ \Psi_\Omega^\delta(u,v) =  \varepsilon_u^2 \dfrac{\Vert\nabla u\Vert^2}{2} + \varepsilon_v^2 \dfrac{\Vert \nabla v\Vert^2}{2} + \int_\Omega F_\delta(u,v)dx. \]
Then, consider \eqref{eq:step3} and pick as test functions $s = \mu_m$, $w = \varphi_m$ while testing the equations for the chemical potentials by $\partial_t u_m := \frac{\partial u_m}{\partial t}$ and $\partial_t v_m := \frac{\partial v_m}{\partial t}$, respectively. From the four resulting equations, owing to the fact that, by definition of orthogonal projection,
		\begin{equation}
			\label{eq:projection}
			(x - \Pi_m(x), w) = 0 \qquad \forall x \in H, \: w \in W_m,
		\end{equation}
		one deduces the energy equality
		\begin{equation} \td{}{t} \Psi_\Omega^\delta(u_m,v_m) + \sigma(v_m-c, \varphi_m)+\|\nabla \mu_m \|^2 + \| \nabla \varphi_m \|^2 = 0.
\label{apprenid}
\end{equation}
		Observe now that
		$$
			\sigma(v_m-c, \varphi_m)   = \sigma\left(v_m-c, \Pi_m\left(-\varepsilon_v^2 \Delta v_m + \pd{F_\delta}{v}(u_m, v_m)\right) \right)
			= \sigma\left(v_m-c, -\varepsilon_v^2 \Delta v_m + \pd{F_\delta}{v}(u_m, v_m) \right),
		$$
		since $v_m-c \in W_m$. Indeed, let us remind that all constant functions on $\Omega$ belong to $W_m$ for any $m \geq 1$, since $w_1 \equiv 1$ is the first eigenfunction of the Neumann--Laplace operator. Therefore, we get
		\begin{align}
\nonumber
		 \sigma(v_m-c, \varphi_m) & = \sigma \varepsilon_v^2 \| \nabla v_m\|^2 + \sigma\left(v_m-c, \pd{F_\delta}{v}(u_m, v_m)\right)\\
			& = \sigma \varepsilon_v^2 \| \nabla v_m\|^2 + \sigma\left(v_m-c, \hat{S}_\delta'(v_m)-\theta_{0,v}v_m + \pd{C}{v}(u_m,v_m)\right).
\label{Oonoterm}		
\end{align}
Observe now that differentiating \eqref{eq:taylor} twice yields
		\begin{equation*}
			\hat{S}_\delta''(r) =
			\begin{cases}
				\displaystyle\sum_{i=0}^{2} \dfrac{\hat{S}^{(i+2)}(-1+\delta)}{i!}[r-(-1+\delta)]^i & \quad \text{if } r \leq -1 + \delta,\\[0.2cm]
				\hat{S}''(r) & \quad \text{if } |r| \leq 1-\delta,\\
				\displaystyle\sum_{i=0}^{2} \dfrac{\hat{S}^{(i+2)}(1-\delta)}{i!}[r-(1-\delta)]^i& \quad \text{if } r \geq 1 - \delta,
			\end{cases}
		\end{equation*}
		therefore $\hat{S}_\delta$ is convex (see properties (iv), (v), (vi)) provided $\delta$ is sufficiently small (i.e. $\delta < \min(\epsilon_0, \epsilon_1)$).
Hence we infer that
		\[ \ii{\Omega}{}{\hat{S}_\delta(v_m)}{x} \leq \ii{\Omega}{}{\hat{S}_\delta'(v_m)(v_m - c)}{x} + \ii{\Omega}{}{\hat{S}_\delta(c)}{x} \]
		or, equivalently,
		\begin{equation}
			\label{eq:p-energy1}
			\ii{\Omega}{}{\hat{S}_\delta'(v_m)(v_m - c)}{x} \geq \ii{\Omega}{}{\hat{S}_\delta(v_m)}{x} - \ii{\Omega}{}{\hat{S}_\delta(c)}{x}.
		\end{equation}
		Moreover, we observe that
		\begin{equation}
			\label{eq:p3}
			Q(u_m, v_m) := (v_m-c)\left(-\theta_{0,v}v_m + \pd{C}{v}(u_m,v_m)\right)
		\end{equation}
		is an algebraic polynomial of degree 3. Thus, on account of \eqref{Oonoterm}-\eqref{eq:p-energy1},  from \eqref{apprenid} we deduce
		\begin{equation}
			\label{eq:energystep2}
			\td{}{t} \Psi_\Omega^\delta(u_m,v_m) + \sigma\varepsilon^2_v\| \nabla v_m \|^2 + \|\nabla \mu_m \|^2 + \| \nabla \varphi_m \|^2 \\+ \sigma \ii{\Omega}{}{Q(u_m,v_m)}{x} \leq \ii{\Omega}{}{\hat{S}_\delta(c)}{x}.
		\end{equation}
		Recalling now Proposition \ref{prop:fdelta}, we can find $K(k_1, |\Omega|, \sigma) > 0$ such that
		\begin{equation}
			\label{eq:energystep5}
			\left| \ii{\Omega}{}{\sigma Q(u_m,v_m)}{x} \right| \leq \ii{\Omega}{}{k_1(u_m^4 + v_m^4)}{x} + K(k_1, |\Omega|, \sigma).
		\end{equation}
		Let us redefine the potential approximation by adding the constant $k_2$ appearing in Proposition \ref{prop:fdelta}, namely
		\[ \widetilde{F}_\delta(u,v) = F_\delta(u,v) + k_2, \]
		so that the very same result entails
		\begin{equation}
			\label{eq:energystep6}
			\widetilde{F}_\delta(u,v) \geq k_1(u^4 + v^4) \qquad \forall \: u,v \in \mathbb{R}.
		\end{equation}
		Accordingly, we set
		\[ \widetilde{\Psi}_\Omega^\delta(u,v) = \dfrac{\varepsilon_u^2}{2} \| \nabla u \|^2 + \dfrac{\varepsilon_v^2}{2} \| \nabla v \|^2 + \int_\Omega \widetilde{F}_\delta(u,v) dx,\]
		in order to have a positive (and coercive) energy functional. Exploiting \eqref{eq:energystep5} and \eqref{eq:energystep6} in \eqref{eq:energystep2}, jointly with
		\[
		\ii{\Omega}{}{\widetilde{F}_\delta(u_m,v_m)}{x}\leq \widetilde{\Psi}_\Omega^\delta(u_m,v_m),
		\]
		we obtain
		\begin{equation}
			\label{eq:energystep7}
			\td{}{t} \widetilde{\Psi}_\Omega^\delta(u_m,v_m) + \sigma\varepsilon^2_v\| \nabla v_m \|^2 + \|\nabla \mu_m \|^2 + \| \nabla \varphi_m \|^2 \\\leq \widetilde{\Psi}_\Omega^\delta(u_m,v_m)+ \ii{\Omega}{}{\hat{S}_\delta(c)}{x} + K.
		\end{equation}
		It is straightforward to see that the constant at right hand side is independent of the Galerkin parameter $m$. An application of the Gronwall lemma yields the energy inequality
		\begin{equation}
			\label{eq:energystep8}
			\widetilde{\Psi}_\Omega^\delta(u_m(t),v_m(t)) \leq \widetilde{\Psi}_\Omega^\delta(u_{0,m},v_{0,m})e^t + K
		\end{equation}
		for some $K > 0$ independent of $m$ (the constant $K$ is redefined without relabeling it). The bound implies that the approximating solutions are well defined on the whole $[0,T]$. Integrating \eqref{eq:energystep7} over $[0,t]$, $t \in (0,T]$, and using \eqref{eq:energystep8}, we eventually get
			\begin{equation}
			\label{eq:finenbound}
			\widetilde{\Psi}_\Omega^\delta(u_m(t),v_m(t)) + \int_0^t \|\nabla \mu_m(\tau) \|^2 + \| \nabla \varphi_m(\tau) \|^2 \: \mathrm{d}\tau \leq K_1\widetilde{\Psi}_\Omega^\delta(u_{0,m},v_{0,m})e^{T} + K_2, \qquad \forall\, t\in [0,T],
		\end{equation}
	for some positive constants $K_1, K_2$ independent of $m$. Existence can now be recovered through \eqref{eq:finenbound}.
Indeed, we can first pass to the limit in the Galerkin scheme to obtain the existence of a solution to the regularized problem (see, e.g. \cite[Sec. 3]{Miranville2}). Then we use uniform estimates with respect $\delta$ to get the existence of a weak solution letting $\delta$ go to $0$ along a suitable sequence (see, for instance, \cite[Sec. 3]{GGW} for the details).

\subsection{Proof of \eqref{eq:dissipative}}
\label{subsec:dissest}
We can now take advantage of the bounds (vii) and (viii) (see Definition \ref{def:solution}) which hold for any given $T>0$.
Indeed, consider \ref{eq:weak2} in Definition \ref{def:solution} and choose $s = \mu$, $w = \varphi$. Then test the equations for the chemical potentials
(see (v)-(vi) in Definition \ref{def:solution}) by $\frac{\mathrm{d}u}{\mathrm{d}t}$ and $\frac{\mathrm{d}v}{\mathrm{d}t}$, respectively. From the four resulting equations, we deduce the energy identity (see, e.g., \cite[Lemma 4.1]{RS})
	\begin{equation}
		\label{eq:energyeq}
		\td{}{t} \Psi_\Omega(u,v) + (\sigma(v-c), \varphi)+\|\nabla \mu \|^2 + \| \nabla \varphi \|^2 = 0.
	\end{equation}
	Observe that
	\begin{equation}
		\label{eq:oono1}
		\begin{split}
		(\sigma(v-c), \varphi) & = (\sigma(v-c), -\varepsilon_v^2\Delta v + \pd{F}{v}(u,v)) =\\ &= (\sigma(v-c), -\varepsilon_v^2\Delta v) + (\sigma(v-c), \hat{S}'(v)) + (\sigma(v-c), - \theta_{0,v} v + \alpha u + 2\beta uv + \gamma u^2).
		\end{split}
	\end{equation}
	Then testing the equation for the chemical potential $\mu$ in \eqref{eq:coupledCH} by $\sigma (u-\overline{u_0})$, we get
	\begin{equation}
		\label{eq:oono2}
		\begin{split}
			(\sigma (u-\overline{u_0}), \mu) & = (\sigma (u-\overline{u_0}), -\varepsilon_u^2\Delta u + \pd{F}{u}(u,v)) =\\ &= (\sigma (u-\overline{u_0}), -\varepsilon_u^2\Delta u) + (\sigma (u-\overline{u_0}), \hat{S}'(u)) + (\sigma (u-\overline{u_0}), -\theta_{0,u} u + \alpha v + \beta v^2 + 2\gamma uv).
		\end{split}
	\end{equation}
	The terms in \eqref{eq:oono1} can be treated as follows (similar considerations hold for \eqref{eq:oono2}).
	\begin{itemize}
		\item[$*$] The first term satisfies
		\[ (\sigma(v-c), -\varepsilon_v^2\Delta v) = \sigma \varepsilon_v^2 \| \nabla v \|^2, \]
		after an integration by parts, and exploiting the boundary conditions for $v$.
		\item[$*$] The second term is treated exploiting the convexity of the singular part $\hat{S}$, namely
		\[ \ii{\Omega}{}{\hat{S}(v)}{x} \leq \ii{\Omega}{}{\hat{S}'(v)(v - c)}{x} + \ii{\Omega}{}{\hat{S}(c)}{x}, \]
		or, equivalently,
		\[  \ii{\Omega}{}{\hat{S}'(v)(v - c)}{x} \geq \ii{\Omega}{}{\hat{S}(v)}{x} - \ii{\Omega}{}{\hat{S}(c)}{x}. \]
		\item[$*$] The third term is treated noticing that, for instance,
		\[ v(-\theta_{0,v} v + \alpha u + 2\beta uv + \gamma u^2) = C(u,v) + \beta uv^2 - \theta_{0,v}v^2.\]
	\end{itemize}
	Collecting all the resulting terms obtained from \eqref{eq:oono1} and \eqref{eq:oono2}, we get	
	\begin{multline*}
		(\sigma(v-c), \varphi) \geq \dfrac{1}{2}\sigma \varepsilon_v^2 \| \nabla v \|^2 + \ii{\Omega}{}{\sigma \hat{S}(v)  - \sigma \dfrac{\theta_{0,v}}{2}v^2 }{x}+ \ii{\Omega}{}{\sigma C(u,v) + \sigma\beta uv^2  - \sigma \dfrac{\theta_{0,v}}{2}v^2 }{x} \\ - \ii{\Omega}{}{\sigma \hat{S}(c) + \sigma c( -\theta_{0,v} v + \alpha u + 2 \beta uv + \gamma u^2)}{x},
	\end{multline*} \vspace{-0.4cm}
	\begin{multline*}
		(\sigma (u-\overline{u_0}), \mu) \geq \dfrac{1}{2} \sigma \varepsilon_u^2 \| \nabla u \|^2 + \ii{\Omega}{}{\sigma \hat{S}(u) - \sigma \dfrac{\theta_{0,u}}{2}u^2}{x} + \ii{\Omega}{}{\sigma C(u,v) + \sigma\gamma u^2v - \sigma \dfrac{\theta_{0,u}}{2}u^2}{x} \\ - \ii{\Omega}{}{\sigma\hat{S}(\overline{u_0}) + \sigma \overline{u_0}( -\theta_{0,u} u + \alpha v + \beta v^2 + 2 \gamma uv)}{x}.
	\end{multline*}
	Let us now turn back to the energy identity \eqref{eq:energyeq} and add to both sides the quantity $(\sigma (u-\overline{u_0}), \mu)$. Exploiting the above inequalities, we find
	\begin{equation*}
		\td{}{t} \Psi_\Omega(u,v) + \sigma \Psi_\Omega(u,v) + \|\nabla \mu \|^2 + \| \nabla \varphi \|^2 \leq (\sigma (u-\overline{u_0}), \mu) + \int_\Omega P(u,v) \: \mathrm{d}x,
	\end{equation*}
	where $P(u,v)$ is a polynomial function of $u$ and $v$. \color{black} With $u$ and $v$ satisfying \color{black} (see Remark \ref{rem:reg3})
\begin{equation}
\label{eq:essbounded}
\esssup{t \in [0,T]}{\Vert u(t)\Vert_{L^\infty}} \leq 1,\qquad  \esssup{t\in [0,T]}{\Vert v(t)\Vert_{L^\infty}} \leq 1,\qquad
\end{equation}
and observing that, thanks to Poincar\'{e}, Young and Cauchy-Schwarz inequalities,
	\[
	(\sigma (u-\overline{u_0}), \mu) = (\sigma (u-\overline{u_0}), \mu - \overline{\mu}) \leq K(1 + \overline{u_0})^2 + \dfrac{1}{2} \| \nabla \mu \|^2,
	\]
	for some positive constant $K$, we get, on any time interval $(0,T)$, the inequality
	\[
		\td{}{t} \Psi_\Omega(u,v) + \sigma \Psi_\Omega(u,v) + \dfrac{1}{2}\|\nabla \mu \|^2 + \| \nabla \varphi \|^2 \leq C,
	\]
	for some positive constant $C$ depending only the parameters of the problem, including the domain and the initial conditions. Then, the estimate follows from the Gronwall lemma.

\subsection{Proof of \eqref{contdep}}
\label{subsec:contdep}

We adapt the argument of \cite{GGM} to the present case. In this subsection $C_i$, $i\in\mathbb{N}$, stands for a positive constant depending only on the parameters of the problem (possibly including $\Omega$ and $T$).
Let us recall that $u_{01}, v_{01}, u_{02}, v_{02} \in V$ are such that $F(u_{01}, v_{01}), F(u_{02}, v_{02})  \in L^1(\Omega)$ and $\overline{u_{01}}, \overline{v_{01}}, \overline{u_{02}}, \overline{v_{02}} \in (-1,1)$. In particular, there exists $m \in (0,1]$ such that
	\begin{equation}
		\label{eq:initialconditions}
		\overline{u_{01}}, \overline{v_{01}}, \overline{u_{02}}, \overline{v_{02}} \in (-m, m).
	\end{equation}
	Furthermore, let $R > 0$ be a positive initial energy bound, i.e.
	\begin{equation}
		\label{eq:initialenergy}
		\max\{\Psi_\Omega(u_{01},v_{01}), \Psi_\Omega(u_{02},v_{02})\}\leq R.
	\end{equation}
	Let us consider the weak formulations of Problem \eqref{eq:coupledCHoc} with initial data $(u_{0i},v_{0i})$, $i=1,2$.
	We know that they admit at least a finite energy solution, say, $(u_i,v_i)$. Set now
	\[
	\begin{split}
		u & = u_1 - u_2, \qquad v = v_1 - v_2,\\
		\mu& = \mu_1 - \mu_2,\qquad \varphi = \varphi_1 - \varphi_2,\\
		u_0 &= u_{01} - u_{02},\quad v_0 = v_{01} - v_{02}.
	\end{split}
	\]
	Then, we formally have
	\begin{equation}
		\label{eq:coupledCHstab2}
		\begin{cases}
			\pd{u}{t} = \Delta \mu & \quad \text{in } \Omega \times (0,T),\\
			\mu = -\varepsilon_u^2\Delta u + \dfrac{\partial  F}{\partial u}(u_1,v_1) - \dfrac{\partial  F}{\partial u}(u_2,v_2)  & \quad \text{in } \Omega \times (0,T),\\
			\pd{v}{t} + \sigma v  = \Delta \varphi & \quad \text{in } \Omega \times (0,T),\\
			\varphi = -\varepsilon_v^2 \Delta v + \dfrac{\partial F}{\partial v}(u_1,v_1) - \dfrac{\partial F}{\partial v}(u_2,v_2)& \quad \text{in } \Omega \times (0,T),\\
			\pd{u}{\mathbf{n}} = \pd{v}{\mathbf{n}} = 0 & \quad \text{on } \partial\Omega \times (0,T),\\[0.2cm]
			\pd{\mu}{\mathbf{n}} = \pd{\varphi}{\mathbf{n}} = 0 & \quad \text{on } \partial\Omega \times (0,T),\\
			u(\cdot, 0) = u_0 & \quad \text{in } \Omega,\\
			v(\cdot, 0) = v_0 & \quad \text{in } \Omega.\\
		\end{cases}
	\end{equation}
	We are now in a position to prove the stability result for Problem \eqref{eq:coupledCHoc}. Firstly, let us point out that the weak formulation of Problem \eqref{eq:coupledCHstab2} reads
	\[
	\begin{cases}
		\left\langle \pd{u}{t}, s \right\rangle + (\nabla \mu, \nabla s) = 0  &\qquad \forall \: s \in V, \: \aevv (0,T),\\[0.3cm]
		\left\langle \pd{v}{t}, w \right\rangle + \sigma(v, w) + (\nabla \varphi, \nabla w) = 0  & \qquad \forall \: w \in V, \: \aevv (0,T).
	\end{cases}	
	\]
	Choosing $s = \mathcal{N}(u - \overline{u})$ and $w = \mathcal{N}(v-\overline{v})$ as test functions, we obtain  	
	\[
	\begin{cases}
		\dfrac{1}{2}\td{}{t}\| u - \overline{u} \|_*^2 + (\mu, u - \overline{u}) = 0 & \qquad \aevv (0,T),\\[0.3cm]
		\dfrac{1}{2}\td{}{t}\| v - \overline{v} \|_*^2 + \sigma \| v - \overline{v} \|_*^2 + (\varphi, v - \overline{v}) = 0 & \qquad \aevv (0,T),
	\end{cases}
	\]
	where we have used Proposition \ref{prop:propertiesAN} and
	\begin{equation}
		\label{eq:propertiesAN}
		\left \langle \pd{u}{t}, \mathcal{N}u \right\rangle = \dfrac{1}{2} \td{}{t}\|u\|_*^2 \: \: \aevv (0,T), \quad \forall \: u \in H^1([0,T], V_0^*).
	\end{equation}
	Next, recalling the definition of the $\| \cdot \|_{-1}$ norm, we find
	\begin{equation}
\begin{cases}
		\label{eq:stab1}
		 		\dfrac{1}{2}\td{}{t}\| u \|_{-1}^2 + (\mu, u - \overline{u}) = 0 , &\qquad \aevv (0,T),\\[0.2cm]
			\dfrac{1}{2}\td{}{t}\| v \|_{-1}^2 + \sigma \| v \|_{-1}^2  + (\varphi, v - \overline{v}) = 0 &\qquad \aevv (0,T),
\end{cases}
    \end{equation}
	on account of the equations for the total masses
	\begin{equation*}
		\label{eq:diffeqavg}
			\td{|\overline{u}|^2}{t} = 0  \quad \aevv (0,T),\qquad
			\dfrac{1}{2}\td{|\overline{v}|^2}{t} + \sigma | \overline{v} |^2 = 0  \quad \aevv (0,T).
	\end{equation*}
The remaining scalar products are managed as follows. For the sake of brevity, we only show how to handle the one appearing in the first equation of \eqref{eq:stab1}. Substituting the expression for $\mu$ (see \eqref{eq:coupledCHstab2}), we get
	\begin{equation*}
		(\mu, u - \overline{u}) = \varepsilon_u^2 \| \nabla u \|^2 + \left(\dfrac{\partial  F}{\partial u}(u_1,v_1) - \dfrac{\partial  F}{\partial u}(u_2,v_2),u\right) - \left(\dfrac{\partial  F}{\partial u}(u_1,v_1) - \dfrac{\partial  F}{\partial u}(u_2,v_2),\overline{u}\right),
	\end{equation*}
	and thus we arrive at (see \eqref{eq:coupling}-\eqref{eq:mixentropy})
	\begin{multline*}
		(\mu, u - \overline{u}) = \varepsilon_u^2 \| \nabla u \|^2 + \left(\hat{S}'(u_1)-\hat{S}'(u_2),u\right) -\left(\hat{S}'(u_1)-\hat{S}'(u_2), \overline{u}\right) \\+ \left(\dfrac{\partial  P_3}{\partial u}(u_1,v_1) - \dfrac{\partial  P_3}{\partial u}(u_2,v_2),u\right) - \left(\dfrac{\partial  P_3}{\partial u}(u_1,v_1) - \dfrac{\partial  P_3}{\partial u}(u_2,v_2),\overline{u}\right),
	\end{multline*}
	where $P_3$ is a third-degree polynomial function such that
	\begin{equation} \label{eq:polybound}
			\begin{split}
			& \left(\dfrac{\partial  P_3}{\partial u}(u_1,v_1) - \dfrac{\partial  P_3}{\partial u}(u_2,v_2),u_1-u_2 \right) \\ & \qquad \qquad = \left( - \theta_{0,u}(u_1-u_2) + \alpha(v_1 - v_2) + \beta(v_1^2 - v_2^2) + 2\gamma(u_1v_1-u_2v_2), u_1 - u_2 \right) \\ & \qquad \qquad \geq - C_1 \|u\|^2 - C_2 \|v\|^2.
			\end{split}
	\end{equation}
	Indeed, we have
	\begin{equation} \label{eq:polynomialproof1}
		(\alpha(v_1 - v_2), u_1 - u_2) \geq - C_3 \| v \|^2 - C_4 \|u\|^2
	\end{equation}
	by the Cauchy-Schwarz and Young inequalities, whereas
	\begin{equation} \label{eq:polynomialproof2}
	\begin{split}
		(\beta(v_1^2 - v_2^2), u_1 - u_2) = (\beta(v_1 - v_2)(v_1+v_2), u_1 - u_2) \geq - C_5 \| v \|^2 - C_6 \|u\|^2,
	\end{split}
	\end{equation}
	and finally
	\begin{equation} \label{eq:polynomialproof3}
	(2\gamma(u_1v_1-u_2v_2), u_1 - u_2) = (2\gamma u_1(v_1-v_2) + 2\gamma v_2(u_1 - u_2), u_1 - u_2) \geq - C_7 \| v \|^2 - C_8 \|u\|^2
	\end{equation}
	where in \eqref{eq:polynomialproof2} and \eqref{eq:polynomialproof3} we have also exploited the essential boundedness of $u_1,v_1,u_2$ and $v_2$ (see \eqref{eq:essbounded}). Bound \eqref{eq:polybound} follows from \eqref{eq:polynomialproof1}-\eqref{eq:polynomialproof3}.
	Moreover, using the essential boundedness once more, we have
	\begin{equation*}
		\left| \left(\dfrac{\partial  P_3}{\partial u}(u_1,v_1) - \dfrac{\partial  P_3}{\partial u}(u_2,v_2),\overline{u}\right) \right| \leq 2|\Omega|(|\alpha| + |\beta| + 2|\gamma|)|\overline{u}|, \qquad\text{ a.e. in } [0,T].
	\end{equation*}
	Finally, from the convexity of $\hat{S}$, it also holds that
	\begin{equation*}
		\left(\hat{S}'(u_1)-\hat{S}'(u_2),u\right) \geq \theta_u \|u\|^2.
	\end{equation*}
	Taking the obtained results into account and arguing similarly for $(\varphi, v-\overline{v})$, from \eqref{eq:stab1} we infer
	\begin{equation}
		\footnotesize
		\label{eq:stab2}
		\begin{cases}
				\dfrac{1}{2}\td{}{t}\| u \|_{-1}^2 + \varepsilon_u^2 \| \nabla u \|^2 + \theta_u \|u\|^2 \leq C_9( \|u\|^2 + \|v\|^2 + |\overline{u}|) + \left(\hat{S}'(u_1)-\hat{S}'(u_2), \overline{u}\right)   & \qquad \aevv (0,T),\\[0.3cm]
				\dfrac{1}{2}\td{}{t}\| v \|_{-1}^2 + \sigma \| v \|_{-1}^2  + \varepsilon_v^2 \| \nabla v \|^2 +\theta_v \|v\|^2 \leq
				C_{10}(\|u\|^2 +  \|v\|^2  + |\overline{v}|) + \left(\hat{S}'(v_1)-\hat{S}'(v_2), \overline{v}\right)& \qquad \aevv (0,T).
		\end{cases}
	\end{equation}
	Thanks to Proposition \ref{prop:propertiesAN}, for any $k,\omega >0$ we can find $C = C(\omega,k) > 0$ such that
	\begin{equation}
		\label{eq:helpaux}
		k \|u\|^2 \leq \omega \| \nabla u \|^2 + C_3 \|u\|^2_{-1}.
	\end{equation}
	In particular, adding together the two inequalities in \eqref{eq:stab2}, multiplying the result by two, and choosing $\omega = \frac{1}{2}\varepsilon_u^2$ (resp. $\frac{1}{2}\varepsilon_v^2$) and using \eqref{eq:helpaux}, we deduce
	\color{black}
	\begin{multline} 		\label{eq:stab3}
			\td{}{t} (\| u \|_{-1}^2 + \| v \|_{-1}^2) + C_{11}(\|u\|^2_{V} + \|v\|^2_{V}) \leq C_{12}( \| u \|_{-1}^2 +  \| v \|_{-1}^2 +
|\overline{u}| + |\overline{v}|)\\ + 2\left(\hat{S}'(u_1)-\hat{S}'(u_2), \overline{u}\right) + 2\left(\hat{S}'(v_1)-\hat{S}'(v_2), \overline{v}\right) \qquad \aevv (0,T).
	\end{multline}	
	\color{black}
	We are only left to deal with the remaining scalar products. Indeed, for instance, we easily notice that
	\[
	\left(\hat{S}'(u_1)-\hat{S}'(u_2), \overline{u}\right) \leq \left( \|\hat{S}'(u_1)\|_{L^1} + \|\hat{S}'(u_2)\|_{L^1}\right)|\overline{u}|,
	\]
	suggesting the necessity of a uniform estimate for the $L^1$-norm of $\hat{S}'$. We recall a well-known property of any singular potential satisfying the properties listed at the beginning of Subsection \ref{subsec:apprF}  namely (see \cite{Miranville3}, see also \cite{FG12})
	\[
	\| \hat{S}'(u_i) \|_{L^1} \leq \mathcal{Q}(|\overline{u_i}|)\left[1+(\hat{S}'(u_i)-\overline{\hat{S}'(u_i)}, u_i - \overline{u_i})\right], \qquad i = 1,2,
	\]
	where $\mathcal{Q}(\cdot)$ is an increasing function. Then, it is straightforward to prove that there exists a constant $C > 0$ such that
	\[ (\hat{S}'(u_i)-\overline{\hat{S}'(u_i)}, u_i - \overline{u_i}) \leq C(\| \nabla \mu_i \| + 1),\qquad i = 1,2,\]
	implying also that
	\[
	\| \hat{S}'(u_i) \|_{L^1} \leq \mathcal{Q}(m)(1+\| \nabla \mu_i \|),  \qquad i = 1,2,
	\]
	for a possibly redefined increasing function $\mathcal{Q}(\cdot)$. Therefore, thanks to the dissipative estimate, the control
	\begin{equation}
		\label{eq:control1}
		\int_0^T \left( \|\hat{S}'(u_1(t))\|_{L^1} + \|\hat{S}'(u_2(t))\|_{L^1}\right) \: \mathrm{d}t \leq C_{13}
	\end{equation}
	holds. Of course, arguing similarly, one also gets
	\[
	\| \hat{S}'(v_i) \|_{L^1} \leq \mathcal{R}(m)(1+\| \nabla \varphi_i \|),
	\]
	where, again, $\mathcal{R}(\cdot)$ is an increasing function, and
	\begin{equation}
		\label{eq:control2}
		\int_0^T \left( \|\hat{S}'(v_1(t))\|_{L^1} + \|\hat{S}'(v_2(t))\|_{L^1}\right) \: \mathrm{d}t  \leq C_{14}.
	\end{equation}
	Thus an application of the Gronwall lemma yields
	\begin{equation*}
		\| u(t) \|_{-1}^2 + \| v(t) \|_{-1}^2 + \int_0^T \| u(t) \|_{V}^2 + \| v(t) \|_{V}^2 \: \mathrm{d}t \leq C_{15}(\| u_0 \|_{-1}^2 + \| v_0 \|_{-1}^2)e^{C_{16}T} + C_{15}(|\overline{u_0}| + |\overline{v_0}|)e^{C_{16}T}.
	\end{equation*}
	Therefore, recalling Proposition \ref{prop:normeq}, we find inequality \eqref{contdep}. We stress that the proof of continuous dependence in the conserved case $c = \overline{v_0}$ is essentially the same, except the third equation in \eqref{eq:coupledCHstab2} reads
	\[
		\pd{v}{t} + \sigma(v-\overline{v_0}) = \Delta \varphi, \quad \text{in } \Omega \times (0,T)
	\]
	and the squared mean satisfies
	\[
	\td{|\overline{v}|^2}{t} = 0,
	\]
	almost everywhere in $(0,T)$.

\section{Regularization properties}
\label{sec:regularity}
Here we consider the conserved case and we show that any weak solution regularizes in finite time. The extension of these results to the off-critical case does not seem straightforward because of the fact that, due to the presence of the coupling term \eqref{eq:coupling}, we cannot use the Galerkin scheme to achieve higher-order estimates on the time derivatives (compare with the single Cahn-Hilliard-Oono equation in \cite{GGM}). This forces us to use difference quotients in time (see, e.g., \cite{GGW}) but this argument does
not work apparently in the off-critical case (see \eqref{offcriticalquot} below).

As we shall see, these regularization effects  are crucial for the investigation of the longtime behavior. The higher-order estimates are obtained working directly on the weak solution instead of using once more the Galerkin scheme (cf. \cite{GGM}). This is due to the fact that the approximated potential $F_\delta$ is not \textit{uniformly} controlled from below because of the coupling term. More precisely, we make use of of the difference quotients in time. More precisely, for any function $f:[0,T] \to X$, $X$ being any real Banach space, for any $h>0$ and any $t\geq 0$, we set
\[\partial_t^h f := \dfrac{f(t+h)-f(t)}{h}. \]

Our first regularity result is given by
\begin{proposition}
	\label{prop:reg1}
	Let the assumptions of Theorem \ref{th:wellposed} hold. Then, for every $\xi > 0$, there exists a constant $K > 0$ depending on the all parameters of the problem such that
	\begin{equation}
\label{eq:timereg}
\|\partial_t u\|_{L^\infty([\xi, t]; V^*)} + \|\partial_t v\|_{L^\infty([\xi, t]; V^*)} + \| \partial_t u \|_{L^2([t,t+1];V)} +\| \partial_t v \|_{L^2([t,t+1];V)} \leq K
\end{equation}
	for every $t \geq \xi$.
\end{proposition}
\begin{proof}
	Observe first that the difference quotients of the weak solutions satisfy the following system of equations
	\begin{equation}
\label{eq:diffquo}
	\begin{cases}
		\langle (\dq{u})_t, s\rangle  + (\nabla \dq{\mu}, \nabla s) = 0 & \quad \forall \: s \in V, \aevv  (0,T),\\
		\langle (\dq{v})_t, w\rangle  + \sigma(\dq{v}, w) + (\nabla \dq{\varphi}, \nabla w) = 0 & \quad \forall \: w \in V, \aevv  (0,T),
	\end{cases}
\end{equation}
	which can be obtained from the weak formulation. The different quotients of the chemical potentials, in turn, enjoy the equations
	\begin{equation} \label{eq:dqpot}
		\begin{cases}
			\dq{\mu} = -\varepsilon_u^2 \Delta \dq{u} + \dfrac{1}{h}\left( \pd{F}{u}(u(t+h),v(t+h))-\pd{F}{u}(u(t),v(t))\right) & \quad \aev,\\
			\dq{\varphi} = -\varepsilon_v^2 \Delta \dq{v} + \dfrac{1}{h}\left( \pd{F}{v}(u(t+h),v(t+h))-\pd{F}{v}(u(t),v(t))\right) & \quad \aev.\\
		\end{cases}
	\end{equation}
	Taking the mass conservation into account, we have
	\[
	\overline{\dq{u}} = \overline{\dq{v}} = 0,
	\]
	or, equivalently,
	\[
	\dq{u}, \dq{v} \in V_0.
	\]

Let us choose as test functions $s = \mathcal{N}\dq{u}$ and $w = \mathcal{N}\dq{v}$. Exploiting Proposition \ref{prop:propertiesAN} and \eqref{eq:propertiesAN}, we get
	\begin{equation}
		\label{eq:reg0}
		\begin{cases}
			\dfrac{1}{2}\td{}{t} \| \dq{u} \|_*^2  + (\dq{\mu}, \dq{u}) = 0 & \quad \aevv (0,T),\\[0.3cm]
			\dfrac{1}{2}\td{}{t} \| \dq{v} \|_*^2  + \sigma\|\dq{v}\|^2_* + (\dq{\varphi}, \dq{v}) = 0 & \quad \aevv (0,T).\\
		\end{cases}
	\end{equation}
Note that in the off-critical case we would have
\begin{equation}
\label{offcriticalquot}
\overline{\dq{v}}= e^{-\sigma t}(\overline{v_0} - c)\frac{e^{-\sigma h} -1}{h}.
\end{equation}
Thus we should take $w = \mathcal{N}(\dq{v}-\overline{\dq{v}})$ and the resulting additional term $(\dq{\varphi}, \overline{\dq{v}})$ seems hard to handle.

	We now deal with the remaining scalar products. Consider, for instance, the scalar product in the first equation of \eqref{eq:reg0}. We have
	\[
	\begin{split}
		(\dq{\mu}, \dq{u}) = \left(\dq{u}, -\varepsilon_u^2 \Delta \dq{u} + \dfrac{1}{h}\left( \pd{F}{u}(u(t+h),v(t+h))-\pd{F}{u}(u(t),v(t))\right) \right),
	\end{split}	
	\]
	and an integration by parts yields (see \eqref{eq:polybound})
	\begin{multline*}
		(\dq{\mu}, \dq{u}) = \varepsilon_u^2 \| \nabla \dq{u} \|^2 + \dfrac{1}{h}\left(\dq{u},\hat{S}'(u(t+h))-\hat{S}'(u(t))\right) \\ + \dfrac{1}{h}\left(\dq{u},\pd{P_3}{u}(u(t+h),v(t+h))-\pd{P_3}{u}(u(t),v(t))\right).
	\end{multline*}
	Note that the second term on the right hand side is non-negative. Furthermore, there holds
	\begin{multline*}
		\pd{P_3}{u}(u(t+h),v(t+h))-\pd{P_3}{u}(u(t),v(t)) = -\theta_{0,u}(u(t+h)-u(t)) + \alpha(v(t+h)-v(t))+ \beta(v^2(t+h)-v^2(t)) \\+ 2\gamma(u(t+h)v(t+h)-u(t)v(t)).
	\end{multline*}
	The nonlinear terms can be treated as in Subsection \ref{subsec:contdep}. In particular, note that
	\[
	\beta(v^2(t+h)-v^2(t)) = \beta(v(t+h)-v(t))(v(t+h)+v(t)),
	\]
	and that, by adding and subtracting $2\gamma u(t)v(t+h)$, we obtain
	\begin{equation*}
		2\gamma(u(t+h)v(t+h)-u(t)v(t)) = 2\gamma v(t+h)(u(t+h)-u(t)) + 2\gamma u(t)(v(t+h)-v(t)).
	\end{equation*}
	Therefore, in light of Remark \ref{rem:reg3} and using repeatedly Young's inequality, we find
	\begin{multline*}
		(\dq{\mu}, \dq{u}) \geq \varepsilon_u^2 \| \nabla \dq{u} \|^2 - \theta_{0,u}\|\dq{u}\|^2 - \dfrac{|\alpha|}{2} \|\dq{u}\|^2 - \dfrac{|\alpha|}{2} \|\dq{v}\|^2 - |\beta|\|\dq{u}\|^2 - |\beta|\|\dq{v}\|^2 \\ - 3|\gamma|\|\dq{u}\|^2 - |\gamma|\|\dq{v}\|^2.
	\end{multline*}
	Thus there exist $K_1, K_2 > 0$ independent of $h$ such that,
	\begin{equation}
		\label{eq:reg1}
		(\dq{\mu}, \dq{u}) \geq \varepsilon_u^2 \| \nabla \dq{u} \|^2 -K_1\|\dq{u}\|^2 - K_2 \|\dq{v}\|^2.
	\end{equation}
	Similarly, we obtain
	\begin{equation}
		\label{eq:reg2}
		(\dq{\varphi}, \dq{v}) \geq \varepsilon_v^2 \| \nabla \dq{v} \|^2 -K_3\|\dq{u}\|^2 - K_4 \|\dq{v}\|^2.
	\end{equation}
	for some positive constants $K_3, K_4$ also independent of $h$. Adding \eqref{eq:reg1} and \eqref{eq:reg2} together, and exploiting \eqref{eq:helpaux} we get
	\begin{equation}
		\label{eq:reg3}
		(\dq{\mu}, \dq{u}) + (\dq{\varphi}, \dq{v}) \geq \dfrac{\varepsilon_u^2}{2}\| \nabla \dq{u} \|^2 +\dfrac{\varepsilon_v^2}{2}\| \nabla \dq{v} \|^2 - K_5 \| \dq{u} \|_*^2 - K_6 \| \dq{v} \|_*^2,
	\end{equation}
	with $K_5, K_6>0$ and independent of $h$. Summing up, from \eqref{eq:reg0}, thanks to \eqref{eq:reg3}, we deduce the inequality
	\begin{equation}
		\label{eq:reg4}
		\dfrac{1}{2}\td{}{t} \| \dq{u} \|_*^2  + \dfrac{1}{2}\td{}{t} \| \dq{v} \|_*^2  + K_7 \| \dq{u} \|_V^2 + K_8 \| \dq{v} \|_V^2  \leq K_5 \| \dq{u} \|_*^2 + |K_6-\sigma| \| \dq{v} \|_*^2,
	\end{equation}
	which holds almost everywhere in $(0,+\infty)$, with $K_7, K_8>0$ independent of $h$. In \eqref{eq:reg4}, the Poincar\'{e}-Wirtinger inequality has also been used. On account of
	\begin{equation} \label{eq:gronwallhy}
		 \| \dq{u} \|_{L^2([t,t+1];V^*)} \leq \| \partial_t u \|_{L^2([t,t+1+h];V^*)}, \qquad
			\| \dq{v} \|_{L^2([t,t+1];V^*)} \leq \| \partial_t v \|_{L^2([t,t+1+h];V^*)},
		 \end{equation}
	the uniform Gronwall lemma and Proposition \ref{prop:normeq} give
	\[ \|\dq{u}\|_{L^\infty([\xi, t]; V^*)} + \|\dq{v}\|_{L^\infty([\xi, t]; V^*)} + \| \dq{u} \|_{L^2([t,t+1];V)} +\| \dq{v} \|_{L^2([t,t+1];V)} \leq K_{9},\]
	for every $t \geq \xi$. The constant $K_{9}$ depends on the parameters of the problem but is independent of $h$ and thus a passage to the limit $h \to 0^+$ entails \eqref{eq:timereg}.
\end{proof}
\noindent
The following result is about spatial regularity of the order parameters (cf. \cite{GGW}).
\begin{proposition}
	\label{prop:reg2}
	Let the assumptions of Theorem \ref{th:wellposed} hold. Then, for every $\xi > 0$, there exists a constant $K > 0$ depending on the all parameters of the problem such that
	\begin{equation}
    \label{eq:spareg}
	\| u \|_{L^\infty([\xi, t]; H^2(\Omega))} + \| v \|_{L^\infty([\xi, t]; H^2(\Omega))} \leq K
	\end{equation}
	for any $t \geq \xi$.
\end{proposition}
\begin{proof}
	Recalling \ref{eq:weak2} in Definition \ref{def:solution}. Then pick $s = \mathcal{N}(\mu - \overline{\mu})$ and $ w = \mathcal{N}(\varphi - \overline{\varphi})$. Owing to Proposition \ref{prop:propertiesAN}, we obtain (derivatives with respect to time are here denoted with $\partial_t$ for simplicity)
	\begin{equation}
		\label{eq:regg}
		\begin{cases}
			\langle \partial_t u, \mathcal{N}(\mu - \overline{\mu}) \rangle + (\mu, \mu-\overline{\mu}) = 0 \\
			\langle \partial_t v, \mathcal{N}(\varphi - \overline{\varphi}) \rangle + \sigma(v-\overline{v_0}, \mathcal{N}(\varphi - \overline{\varphi}))+  (\varphi, \varphi-\overline{\varphi}) = 0
		\end{cases}
	\end{equation}
	almost everywhere in $(0,+\infty)$. Observe that
	\[
		(\mu, \mu-\overline{\mu}) = \| \mu-\overline{\mu} \|^2,\qquad (\varphi, \varphi-\overline{\varphi}) = \| \varphi-\overline{\varphi} \|^2,
	\]
	since $(c, \mu-\overline{\mu}) = (c, \varphi-\overline{\varphi}) = 0$ for any $c \in \mathbb{R}$. Furthermore, we have
	\begin{equation}
		\label{eq:reg21}
		|\langle \partial_t u, \mathcal{N}(\mu - \overline{\mu}) \rangle| \leq \|\partial_t u\|_{*}\| \mathcal{N}(\mu-\overline{\mu})\|_{V_0},
	\end{equation}
	with $\| v \|_{V_0}^2 := \langle Av, v \rangle$ for any $v \in V_0$. Hence we get
	\begin{equation*}
		\begin{split}
			|\langle \partial_t u, \mathcal{N}(\mu - \overline{\mu}) \rangle| & \leq \|\partial_t u\|_{*}\| \mathcal{N}(\mu-\overline{\mu})\|_{V_0} \\
			& \leq K_1\|\partial_t u\|_{V^*}\| \mu-\overline{\mu}\|\\
			& \leq \dfrac{1}{2}\|\mu-\overline{\mu}\|^2 + K_2\|\partial_t u\|_{V^*}^2,\\
		\end{split}	
	\end{equation*}
	for some positive constants $K_1$ and $K_2$. Here Young's inequality has been used. Arguing similarly, we find
	\begin{equation*}
		\begin{split}
			|\langle \partial_t v, \mathcal{N}(\varphi - \overline{\varphi}) \rangle| & \leq \|\partial_t v\|_{*}\| \mathcal{N}(\varphi -\overline{\varphi })\|_{V_0} \\
			& \leq K_1\|\partial_t v\|_{V^*}\| \varphi -\overline{\varphi}\|\\
			& \leq \dfrac{1}{4}\|\varphi -\overline{\varphi }\|^2 + K_3\|\partial_t v\|_{V^*}^2.\\
		\end{split}	
	\end{equation*}
	Also, we have
	\begin{equation*}
		\begin{split}
			\sigma|(v-\overline{v_0}, \mathcal{N}(\varphi - \overline{\varphi}))| \leq K_4\|v-\overline{v_0}\|^2 + \dfrac{1}{4}\|\varphi - \overline{\varphi}\|^2
		\end{split}
	\end{equation*}
	for some $K_4>0$. Adding the two equations in \eqref{eq:regg} together and using the above inequalities, we find
	\begin{equation}
		\dfrac{1}{2}\| \mu - \overline{\mu} \|^2 + \dfrac{1}{2}\| \varphi - \overline{\varphi} \|^2 \leq K_2\|\partial_t u\|_{V^*}^2+K_3\|\partial_t v\|_{V^*}^2 + K_4\|v-\overline{v_0}\|^2,
	\end{equation}
	which entails (see  \eqref{eq:ENERID}, \eqref{eq:dissipative}, and Proposition \ref{prop:reg1})
    \[
	\| \mu - \overline{\mu} \|_{L^\infty([\xi, t];H)} + \| \varphi - \overline{\varphi} \|_{L^\infty([\xi, t];H)}\leq K_5 \quad \forall \: t \geq \xi,
	\]
	for some positive constant $K_5$. Consider now the nonlinear Neumann problem
	\begin{equation}
		\label{eq:nnlp}
		\begin{cases}
			-\rho\Delta g + \hat{S}'(g) = f & \quad \text{in } \Omega,\\
			\pd{g}{\mathbf{n}} = 0 & \quad \text{on } \partial\Omega,\\
		\end{cases}
	\end{equation}
where $f = \mu - \pdd{P_3}{u}(u,v), \rho = \varepsilon_u^2$ or $f = \varphi - \pdd{P_3}{v}(u,v), \rho = \varepsilon_v^2$. Then, recalling \cite[Lemma 7.1]{GGW} and Remark \ref{rem:reg3}, \eqref{eq:spareg} follows.
\end{proof}
\noindent
Taking advantage of the previous results, a further regularity can be proven (cf. \cite{GGW}), namely,
\begin{proposition}
	\label{prop:reg3}
	Let the assumptions of Theorem \ref{th:wellposed} hold. Let $1 \leq p \leq 6$ if $d = 3$, and $1 \leq p < +\infty$ if $d = 2$. Then, for every $\xi > 0$, there exists a constant $K > 0$ depending on the all parameters of the problem (including $p$) such that
	\[
	\| \mu \|_{L^\infty([\xi, t]; V)} + \| \varphi \|_{L^\infty([\xi, t]; V)} + \| u \|_{L^\infty([\xi, t]; W^{2,p}(\Omega))} + \| v \|_{L^\infty([\xi, t]; W^{2,p}(\Omega))} \leq K
	\]
	for every $t \geq \xi$. Furthermore, there exists $L > 0$, also depending on the all parameters of the problem (including $p$), such that
	\[ \left\|\pd{F}{u}(u,v)\right\|_{L^\infty([\xi, t]; L^p(\Omega))} + \left\|\pd{F}{v}(u,v)\right\|_{L^\infty([\xi, t]; L^p(\Omega))} \leq L\]
	for every $t \geq \xi$.
\end{proposition}
\begin{proof}
In this proof $K$ denotes a generic positive constant depending at most on the all parameters of the problem (including $p$).
	Let us recall that (see Subsection \ref{contdep})
	\[	
		\| \hat{S}'(u) \|_{L^1} \leq K(1+\| \nabla \mu \|), \qquad
		\| \hat{S}'(v) \|_{L^1} \leq K(1+\| \nabla \varphi \|).
	\]
	Observe that
	\[ \| \mu \|_V \leq \| \mu - \overline{\mu} \|_V + \| \overline{\mu} \|_V = \| \mu - \overline{\mu} \|_V + |\Omega|^\frac{1}{2}|\overline{\mu}|. \]
	\color{black}Moreover\color{black}, testing the equation for $\mu$ with the characteristic function $\chi_\Omega$, we get (see Remark \ref{rem:reg3})
	\[
		|\Omega| |\overline{\mu}|  = \left| \left(\pd{F}{u}(u,v),1 \right) \right|
		\leq \left| \left(\hat{S}'(u),1 \right) \right|+ \left| \left(\pd{P_3}{u}(u,v),1 \right) \right|
		\leq \| \hat{S}'(u) \|_{L^1} + K.
	\]
	\color{black}
	From the previously mentioned inequalities, we then deduce
	\[
	|\overline{\mu}| \leq K(1 + \| \nabla \mu \|), \qquad 	|\overline{\varphi}| \leq K(1 + \| \nabla \varphi \|).
	\]
	Fix any $\xi > 0$ and $t \geq \xi$. Let us pick $s = \mu$ and $w = \varphi$ in the weak formulation given in Definition \ref{def:solution}-(ix), getting on one hand
	\[
	\|\nabla \mu\|^2 = - \langle \partial_t u, \mu \rangle, \]
	Making use of Proposition \ref{prop:reg1} and applying Young's inequality we have that
	\[
	\|\nabla \mu\|^2 =|\langle \partial_t u, \mu \rangle| \leq \|\partial_t u\|_{V^*}\|\mu\|_V \leq C(1+\|\nabla \mu\|) \leq C + \dfrac{1}{2}\|\nabla \mu\|^2
	\]
	almost everywhere in $[\xi, t]$. Therefore $\nabla \mu \in L^\infty([\xi, t]; H)$.\\
	On the other hand, a similar argument, still thanks to Proposition \ref{prop:reg1} jointly with the Young inequality, yields
	\[
	\|\nabla \varphi\|^2 = -\langle \partial_t v, \varphi \rangle - \sigma(v-\overline{v_0}, \varphi)\leq \|\partial_t v\|_{V^*}\|\varphi\|_V + C\|v-\overline{v_0}\|\|\varphi\|_V\leq C(1+\|\nabla \varphi\|) \leq C + \dfrac{1}{2}\|\nabla \varphi\|^2,
	\]
	\color{black}
	and thus we conclude that
	\[
	\| \mu \|_{L^\infty([\xi, t]; V)} + \| \varphi \|_{L^\infty([\xi, t]; V)} \leq K.
	\]
	Recalling \eqref{eq:nnlp}, thanks to \cite[Lemma 7.4]{GGW} and Remark \ref{rem:reg3}, we also learn that
	\[
	\| u \|_{L^\infty([\xi, t]; W^{2,p}(\Omega))} + \| v \|_{L^\infty([\xi, t]; W^{2,p}(\Omega))} \leq K,
	\]
	and
	\[
	\| \hat{S}'(u)\|_{L^\infty([\xi, t]; L^p(\Omega))} + \| \hat{S}'(v) \|_{L^\infty([\xi, t]; L^p(\Omega))} \leq K.
	\]
    We also know that (see 	Remark \ref{rem:reg3})
	\[ \left\|\pd{P_3}{u}(u,v)\right\|_{L^\infty([\xi, t]; L^p(\Omega))} + \left\|\pd{P_3}{v}(u,v)\right\|_{L^\infty([\xi, t]; L^p(\Omega))} \leq K.
	\]
	This concludes the proof.
\end{proof} \noindent

\begin{remark}
	It is easy to check that indeed, for any $\xi > 0$,
	\[
	\mu(t), \varphi(t) \in H^2(\Omega), \qquad \forall \: t \geq \xi.
	\]
	Thus, the weak solutions are indeed strong, i.e., they satisfy the equations of \eqref{eq:coupledCH} almost everywhere in $\Omega \times (\xi,+\infty)$.
\end{remark}

	\section{The strict separation property in two dimensions} \label{sec:separation}
	In the investigation of Cahn-Hilliard type equations, a very interesting issue concerns the so-called strict separation property, that is, the order parameter stays uniformly away from the pure phases $\pm1$. In other words, entropy always prevails to a certain degree. The validity of this property for Cahn-Hilliard equations with constant mobility has been proven so far in the two-dimensional case only (see \cite{Miranville3}, see also \cite{GGM}).
Here we show that, in two dimensions, the strict separation property holds both for the macrophase and the microphase. More precisely, for any time $t_0>0$, there exist $\omega_u, \omega_v \in (0,1)$ depending on $t_0$, such that, for all $t\geq t_0$,
	\begin{equation*}
		 \|u(t)\|_{L^\infty} \leq 1-\omega_u  ,\qquad \|v(t)\|_{L^\infty} \leq 1-\omega_v
		 \end{equation*}
Our basic assumption is (cf. \cite{GGM})
	\begin{assumption}
		The singular part of the potential $\hat{S}$ satisfies
		\[
			|\hat{S}''(r)| \leq e^{C|\hat{S}'(r)|+C}, \quad \forall \: r \in (-1,1)
		\]
		and is such that $\hat{S}''$ is convex.
	\end{assumption}
	\begin{remark}
		The logarithmic potential \eqref{eq:singularpart} satisfies Assumption A.
	\end{remark} \noindent
We establish the strict separation property by adapting a method developed in \cite{GGM} (see also \cite{GGW,Miranville3}). Let us firstly prove some preliminary results.
	\begin{lemma} \label{lem:sep1}
		Suppose $d=2$. Let the hypotheses of Theorem \ref{th:wellposed} and Assumption A hold. Then, for every $1 \leq p < +\infty$ there exists a positive constant $C$ (depending also on $p$) such that
		\[
		\| \hat{S}''(u) \|_{L^\infty(\xi, t; L^p(\Omega))} + \| \hat{S}''(v) \|_{L^\infty(\xi, t; L^p(\Omega))} \leq C,
		\]
		for every $t \geq \xi$.
	\end{lemma}
	\begin{proof}
		Consider the semilinear Neumann problem \eqref{eq:nnlp}. Choosing, as performed in the proof of Proposition \ref{prop:reg3}, $f = \mu - \pdd{P_3}{u}$ and $\rho = \varepsilon_u^2$, then, owing to Assumption A, \cite[Lemma 7.4]{GGW} entails
		\[ \| \hat{S}''(u) \|_{L^p} \leq C(1 + e^{C\|\mu - \pdd{P_3}{u}\|^2_V}). \]
		Thanks to  Remark \ref{rem:reg3} and Proposition \ref{prop:reg3}, observing that
		\[ \nabla \pd{P_3}{u} = \alpha \nabla v + 2\beta v \nabla v + 2\gamma u \nabla v + 2\gamma v \nabla u,\]
		on account of \eqref{eq:ENERID} and \eqref{eq:dissipative}, it is straightforward to conclude that there exists a constant $C_1$, depending on all the parameters of the problem, $p$ and $\xi > 0$, such that
		\[ \| \hat{S}''(u(t)) \|_{L^p} \leq C_1 \quad \: \forall \: t \geq \xi.\]
		Hence, we get
		\[ \| \hat{S}''(u) \|_{L^\infty(\xi, t; L^p(\Omega))} \leq C_1. \]
		A similar bound holds for $\hat{S}''(v)$, choosing $f = \varphi - \pdd{P_3}{v}$ and $\rho = \varepsilon_v^2$.
	\end{proof} \noindent
	We can now gain higher regularity for the temporal derivatives of $u$ and $v$ as well as for the chemical potentials. Indeed we have
	\begin{lemma} \label{lem:sep2}
		Suppose $d=2$. Let the hypotheses of Theorem \ref{th:wellposed} and Assumption A hold. Then, for any $\xi > 0$, there exists a positive constant $C$, depending on all the parameters of the problem and $\xi$, such that
		\[
		\|\partial_t u \|_{L^\infty(2\xi, t; H)} + \|\partial_t v \|_{L^\infty(2\xi, t; H)} + \|\mu \|_{L^\infty(2\xi, t; H^2(\Omega))} + \|\varphi \|_{L^\infty(2\xi, t; H^2(\Omega))} \leq C,
		\]
		for every $t \geq 2\xi$.
	\end{lemma}
	\begin{proof}
		We work with difference quotients as in the proof of Proposition \ref{prop:reg1}. In this proof $C$ will stand for a generic positive constant depending on the parameters of the problem but independent of $t$ and $h>0$.  Choosing $s = \dq{u}$ and $w = \dq{v}$ in \eqref{eq:diffquo} yields
	\begin{equation}
		\label{eq:sep0}
		\begin{cases}
			\dfrac{1}{2}\td{}{t} \| \dq{u} \|^2 + (\nabla \dq{\mu}, \nabla \dq{u}) = 0 & \quad \aev,\\[0.3cm]
			\dfrac{1}{2}\td{}{t} \| \dq{v} \|^2  + \sigma\|\dq{v}\|^2 + (\nabla \dq{\varphi}, \nabla \dq{v}) = 0 & \quad \aev.\\
		\end{cases}
	\end{equation}
	Integrating by parts the first equation we obtain (see \eqref{eq:dqpot})
	\begin{multline*}
		(\nabla \dq{\mu}, \nabla \dq{u}) = - (\dq{\mu}, \Delta \dq{u}) = \varepsilon_u^2 \| \Delta \dq{u} \|^2 - \dfrac{1}{h}\left(\hat{S}'(u(t+h)) - \hat{S}'(u(t)), \Delta \dq{u} \right) \\ - \dfrac{1}{h}\left(\pd{P_3}{u}(u(t+h), v(t+h)) - \pd{P_3}{u}(u(t), v(t)), \Delta \dq{u} \right).
	\end{multline*}
	Arguing as in Proposition \ref{prop:reg1}, we deduce that
	\[
	\left |\dfrac{1}{h}\left(\pd{P_3}{u}(u(t+h), v(t+h)) - \pd{P_3}{u}(u(t), v(t)), \Delta \dq{u} \right) \right | \leq C \left(\| \nabla \dq{u} \|^2 + \| \nabla \dq{v} \|^2\right)
	\]
	The singular term is treated exploiting the convexity of $\hat{S}''$, owing to Assumption A. In particular, arguing as in \cite[Lemma 5.2]{GGM}, we can prove that
	\[
	\left| \dfrac{1}{h}\left(\hat{S}'(u(t+h)) - \hat{S}'(u(t)), \Delta \dq{u} \right) \right| \leq C \left( \|\hat{S}''(u(t+h))\|_{L^3}^2 + \|\hat{S}''(u(t))\|_{L^3}^2\right) \| \dq{u}\|^2_{L^6} + \dfrac{\varepsilon_u^2}{2} \| \Delta \dq{u} \|^2.
	\]
Thanks to the Sobolev embedding $H^1(\Omega) \hookrightarrow L^6(\Omega)$ and to the Poincar\'{e} inequality it is possible to infer that
	\[ \| \dq{u}\|^2_{L^6} \leq C \|\nabla \dq{u} \|^2 \leq C \| \dq{u} \| \| \Delta \dq{u} \|, \]
	where the last inequality follows from the Cauchy--Schwarz inequality after an integration by parts. Collecting the results, owing also to Young's inequality, we end up with
	\begin{equation} \label{eq:sep1}
		\dfrac{1}{2}\td{}{t} \| \dq{u} \|^2 + \dfrac{1}{4}\| \Delta \dq{u} \|^2 \leq C\left[ \left(\|\hat{S}''(u(t+h))\|_{L^3}^4 + \|\hat{S}''(u(t))\|_{L^3}^4\right)\|\dq{u}\|^2 +  \| \nabla \dq{u} \|^2 + \| \nabla \dq{v} \|^2\right],
	\end{equation}
	almost everywhere in $(0,+\infty)$. Arguing similarly, we deduce
	\begin{multline} \label{eq:sep2}
		\dfrac{1}{2}\td{}{t} \| \dq{v} \|^2 + \sigma \| \dq{v} \|^2 + \dfrac{1}{4}\| \Delta \dq{v} \|^2 \leq C\left(\|\hat{S}''(v(t+h))\|_{L^3}^4 + \|\hat{S}''(v(t))\|_{L^3}^4\right)\|\dq{v}\|^2 \\ +  C \left(\| \nabla \dq{u} \|^2 +  \| \nabla \dq{v} \|^2 \right),
	\end{multline}
	almost everywhere in $(0,+\infty)$. Setting now
	\begin{gather*}
		Y(t) := \|\dq{u}\|^2 + \|\dq{v}\|^2, \\
		W(t) :=  C\left(1 + \|\hat{S}''(u(t+h))\|_{L^3}^4 + \|\hat{S}''(u(t))\|_{L^3}^4 + \|\hat{S}''(v(t+h))\|_{L^3}^4 + \|\hat{S}''(v(t))\|_{L^3}^4\right),
	\end{gather*}
and adding \eqref{eq:sep1} and \eqref{eq:sep2} together, we get the differential inequality (see also \eqref{eq:ENERID} and \eqref{eq:dissipative})
	\begin{equation}
		\dfrac{1}{2} \td{}{t} Y(t) + \dfrac{1}{4}\| \Delta \dq{u} \|^2 + \dfrac{1}{4}\| \Delta \dq{v} \|^2 \leq W(t)Y(t), \qquad \aev.
	\end{equation}
	By virtue of Lemma \ref{lem:sep1}, we have
	\begin{equation}
		\int_t^{t+1} W(\tau) \: \mathrm{d}\tau \leq C, \qquad \forall \: t \geq \xi.
	\end{equation}
	Thus, recalling \eqref{eq:gronwallhy}, an application of the uniform Gronwall lemma and passage to the limit as $h \to 0^+$ entail that, for any $\xi > 0$,
	\[
	\|\partial_t u \|_{L^\infty(2\xi, t; H)} + \|\partial_t v \|_{L^\infty(2\xi, t; H)} \leq C,
	\]
	for every $t \geq \xi$. A comparison argument in the evolution equations of Problem \eqref{eq:coupledCH} yields
	\[\|\Delta \mu\|_{L^\infty(2\xi, t; H)} + \| \Delta \varphi \|_{L^\infty(2\xi, t; H)} \leq C.\]
	The elliptic regularity theory and Proposition \ref{prop:reg3} allow us to get the wanted bound.
	\end{proof} \noindent
	We can now prove the strict separation property for both $u$ and $v$.
	\begin{proposition}
		Suppose $d=2$. Let the hypotheses of Theorem \ref{th:wellposed} and Assumption A hold. Then, for every $\xi > 0$ there exist $\omega_u, \omega_v \in(0,1)$ such that
		\begin{equation*}
				\|u(t)\|_{L^\infty} \leq 1-\omega_u ,\qquad \|v(t)\|_{L^\infty} \leq 1-\omega_v,
		\end{equation*}
	for all $t\geq 2\xi$.
	\end{proposition}
	\begin{proof}
		Arguing on the nonlinear Neumann problem \eqref{eq:nnlp} as in \cite[Lemma 7.2]{GGW} we infer that
		\[
		\| \hat{S}'(u) \|_{L^\infty}+ \|\hat{S}'(v)\|_{L^\infty} \leq C_1\left(\| \mu \|_{L^\infty} + \| \varphi \|_{L^\infty} + \left\| \pd{P_3}{u} (u,v) \right\|_{L^\infty} + \left\| \pd{P_3}{v} (u,v) \right\|_{L^\infty}\right).
		\]
		Owing to Lemma \ref{lem:sep2}, the right hand side is uniformly bounded. Thus the strict separation follows from the properties of $\hat{S}'$.
	\end{proof}
	\begin{remark}
\label{2REG}
		It is now straightforward to prove that $u,v$ also belong to $L^\infty(2\xi, t; H^4(\Omega))$ (see \cite[Cor. 5.1]{GGM}). Thus, via a bootstrap method, provided the boundary of $\Omega$ and the potential $F$ are smooth, we can prove that $u$ and $v$ are as smooth as we want  (see \cite[Rem. 5.2]{GGM}).
	\end{remark} \noindent
	\section{Longtime behavior} \label{sec:equilibrium}
	The main result of this section is the convergence of any finite energy weak solution to a single equilibrium. We will adapt the method exploited in \cite{AW07} for a single Cahn-Hilliard equation. We shall need the regularization properties of weak solutions. Therefore we suppose $c = \overline{v_0}$ (i.e., conserved case). Moreover, without loss of generality, we can take $\overline{u}_0 = \overline{v_0} = 0$ (see Remark \ref{rem:nonzeromean}). Let us restate our (formal) problem in the following equivalent form
	\begin{equation} \label{eq:coupledCHequi}
			\begin{cases}
				\pd{u}{t} = \Delta \mu & \quad \text{in } \Omega \times (0,T),\\
				\mu = - \varepsilon_u^2 \Delta u + \pd{F}{u}(u,v) & \quad \text{in } \Omega \times (0,T),\\
				\pd{v}{t} = \Delta \widetilde{\varphi} & \quad \text{in } \Omega \times (0,T),\\
				\widetilde{\varphi} = - \varepsilon_v^2 \Delta v + \pd{F}{v}(u,v) + \sigma \mathcal{N}v & \quad \text{in }  \Omega \times (0,T),\\
				\pd{u}{\mathbf{n}} = \pd{v}{\mathbf{n}} = 0 & \quad \text{on } \partial \Omega \times (0,T),\\[0.2cm]
				\pd{\mu}{\mathbf{n}} = \pd{\varphi}{\mathbf{n}} = 0 & \quad \text{on } \partial \Omega \times (0,T),\\[0.2cm]
				u(\cdot, 0) = u_0 & \quad \text{in } \Omega,\\
                v(\cdot, 0) = v_0 & \quad \text{in } \Omega.
			\end{cases}
	\end{equation}
	\color{black}Thanks to the conservation of mass, the problem can be viewed as the gradient flow generated by the Ohta-Kawasaki functional
\eqref{eq:ohtakawasaki}. Thus the phase $v$ satisfies a Cahn-Hilliard system in which the chemical potential $\widetilde{\varphi}$ incorporates the reaction term as a nonlocal term.  The resulting energy functional associated with \eqref{eq:coupledCHequi} is
	\[ \ene(u,v) = \Psi_\Omega(u,v) + \dfrac{\sigma}{2}\|v - \overline{v}\|_*^2.\]
	\color{black} Therefore, recalling \eqref{eq:energyeq} and the fact that $c = \overline{v_0} = 0$, we now have the energy identity
	\begin{equation} \label{eq:energyeq2}
		\td{}{t}\ene(u,v) + \|\nabla \mu\|^2 + \|\nabla \widetilde{\varphi}\|^2 = 0
	\end{equation}
	which clearly shows the dissipative nature of the above problem and it is very helpful to investigate the longtime behavior of its solutions.
	\color{black}
	Let us introduce the Hilbert triplet
	\[ V_0 \hookrightarrow H_0 \cong H_0^* \hookrightarrow V_0^*, \]
	where $H_0$ denotes the subspace of $L^2(\Omega)$ functions with null spatial average.
The notion of equilibrium or stationary solution is given by
	\begin{definition} \label{def:statsol}
		A pair $(u_\infty, v_\infty) \in H^2(\Omega)^2$ is a stationary solution to \eqref{eq:coupledCHequi} if
		\begin{equation} \label{eq:coupledCHstat}
			\begin{cases}
				- \varepsilon_u^2 \Delta u_\infty + \pd{F}{u}(u_\infty,v_\infty) = \mu_\infty, & \quad \text{in } \Omega,\\[0.2cm]
				- \varepsilon_v^2 \Delta v_\infty + \pd{F}{v}(u_\infty,v_\infty) + \sigma \mathcal{N}(v_\infty-\overline{v_0}) = \widetilde{\varphi}_\infty & \quad \text{in }  \Omega ,\\
				\pd{u_\infty}{\mathbf{n}} = \pd{v_\infty}{\mathbf{n}} = 0 & \quad \text{on } \partial \Omega ,\\[0.2cm]
				\overline{u_\infty} = \overline{u_0}, \qquad \overline{v_\infty} = \overline{v_0},
			\end{cases}
		\end{equation}
		where $\mu_\infty, \widetilde{\varphi}_\infty \in \mathbb{R}$.
	\end{definition} \noindent

	In the following, all Banach spaces $X^2$, where $X$ is a real Banach space, are considered to be normed with the standard Euclidean norm. Let us consider the set
	\[ Z := \{ (u,v) \in V^2_0\,:\,  \ene(u,v) < +\infty \},\]
	and define the operators $\mathcal{S}(t): Z \to Z$, acting as follows
	\[ \mathcal{S}(t)(u_0,v_0) = (u(t;u_0), v(t;v_0)) \equiv (u(t), v(t)), \qquad t \geq 0. \]
    Note that $Z$ is a complete metric space with respect to the metric induced by the norm in $V_0^2$.
	
    First of all, let us state the following
	\begin{lemma} \label{lem:identity}
		The energy identity
		\begin{equation}
\label{eq:eneriden}
		\ene(u(t),v(t)) + \int_0^t \left(\| \nabla \mu(\tau) \|^2 + \| \nabla \widetilde{\varphi}(\tau) \|^2 \right) \mathrm{d}\tau = \ene(u_0,v_0)
		\end{equation}
		holds for any $t \geq 0$.
	\end{lemma}
	\begin{proof}
	\color{black}	We just need to integrate \eqref{eq:energyeq2} with respect to time over $[0,t]$. \color{black}
	\end{proof} \noindent
	Lemma \ref{lem:identity} entails that $\mathcal{S}(t)(u_0,v_0) \in Z$ for any $t \geq 0$. Moreover, we have $\mathcal{S}(\cdot)(u_0,v_0)  \in \mathcal{C}^0([0,+\infty); Z)$. Also, we can prove that $\mathcal{S}(t) \in \mathcal{C}^0(Z;Z)$ \color{black}(see \cite[Prop. 6.1]{GGM})\color{black}. Thus $(Z,\mathcal{S}(t))$ is a (dissipative) dynamical system. Given $(u_0,v_0) \in Z$, we define the $\omega$-limit set $\omega(u_0,v_0)$ as
	\[
		\omega(u_0,v_0) := \{ (u_\infty, v_\infty) \in H^{2r}(\Omega)^2 \cap Z : \exists \{t_n\}_{n \in \mathbb{N}} \nearrow +\infty \text{ such that } (u(t_n), v(t_n)) \to (u_\infty, v_\infty) \text{ in } H^{2r}(\Omega)^2 \},
	\]
	for $r \in [1/2,1)$. For any $\xi > 0$, Proposition $\ref{prop:reg2}$ entails that the orbits $\{(u(t), v(t))\}_{t \geq \xi}$ are relatively compact in $H^{2r}(\Omega)^2$ for $r \in [1/2,1)$. Thus, we conclude that $\omega(u_0,v_0)$ is a non-empty connected subset of $Z$, and furthermore, by definition, $\omega(u_0,v_0)$ is a compact subset of $H^{2r}(\Omega)^2$ for every $r \in [1/2,1)$ and
	\[ \dist(\mathcal{S}(t)(u_0,v_0), \omega(u_0,v_0)) \to 0 \quad \text{as }t \to +\infty \]
	in the $H^{2r}$-sense (see \cite{AW07}). Next, we prove a second preliminary result.
	\begin{lemma} \label{lem:equi2}
		The functional $\ene(u,v)$ is a strict Lyapunov functional for $\mathcal{S}$, namely energy is conserved only along constant trajectories.
	\end{lemma}
	\begin{proof}
		Lemma \ref{lem:identity} implies that if $\ene(\mathcal{S}(t)(u_0,v_0)) = \ene(u_0,v_0)$ for all $t > 0$, then  $\nabla \mu(t) = \nabla \varphi(t) \equiv 0$ for every $t>0$. Thus $\mathcal{S}(t)(u_0,v_0) = (u_0,v_0)$ for every $t \geq 0$.
	\end{proof} \noindent
	Let $\mathcal{E}$ denote the set of stationary points of $\mathcal{S}$, namely
	\[ \mathcal{E}:= \{(u_0,v_0) \in Z\,:\, \mathcal{S}(t)(u_0,v_0) = (u_0,v_0) \text{ for all } t \geq 0 \}.\]
	As a consequence of Lemma \ref{lem:equi2}, we can prove that $\omega(u_0,v_0) \subset \mathcal{E}$ (see \cite[Theorem 9.2.7]{CH98}). Following \cite{AW07}, we can further characterize the set of stationary points $\mathcal{E}$ as the set of stationary solutions, namely
	\[ \mathcal{E} = \{ (u,v) \in H^2(\Omega)^2 \cap Z\,:\, (u,v) \text{ is a stationary solution}\}.\]
	We now state a (strict) separation property for stationary solutions which follows by adapting the proof of \cite[Proposition 6.1]{AW07}.
	\begin{proposition} \label{prop:maxprinc}
		Let $r \in (d/4, 1)$. For every $f = (f_1, f_2) \in \mathcal{E}$ there exist two constants $M_1, M_2$ such that
		\[ -1 < M_1 \leq f_1(x), f_2(x) \leq M_2 < 1,\]
		for all $x \in \overline{\Omega}$. Furthermore, there exist two constants $K_1, K_2$, independent of $f$, such that
		\[ -1 < K_1 \leq f_1(x), f_2(x) \leq K_2 < 1,\]
		for all $f = (f_1, f_2) \in \omega(u_0,v_0)$ and $x \in \overline{\Omega}$.
	\end{proposition}
	We recall that, by compactness of the $\omega$-limit set in $H^{2r}(\Omega)$, with $r \in (d/4, 1)$, there exists an open set $U_1$ covering $\omega(u_0, v_0)$ such that
	\[ -1 < K_1 - \epsilon < f_1(x), f_2(x) < K_2 + \epsilon < 1, \]
	for every $f = (f_1, f_2) \in U_1$ and some $\epsilon > 0$ independent of $f$. Moreover, $U_1$ attracts the trajectories of the system, since $\omega(u_0,v_0)$ does. Along the lines of \cite{AW07}, it is possible to redefine the free energy $F$. Indeed, let $Q_\epsilon := [K_1-\epsilon, K_2 + \epsilon]^2$ and set
	\[ F_{\text{reg}}(s_1,s_2) = F(s_1,s_2)\chi_{Q_\epsilon}(s_1,s_2) + G(s_1,s_2)\chi_{Q_\epsilon^C}(s_1,s_2),\]
	where $G(s_1,s_2)$ is chosen in such a way to extend $F$ outside $Q_\epsilon$ with $\mathcal{C}^3(\mathbb{R}^2)$ regularity and bounded derivatives up to order three. Accordingly, we introduce the regularized energy functional $\eneg: V_0^2 \to \mathbb{R}$
	\begin{equation}
		\label{eq:ohtakawasaki2}
		\eneg(u, v) =  \varepsilon_u^2 \dfrac{\Vert \nabla u\Vert^2}{2} + \varepsilon_v^2 \dfrac{\Vert \nabla v\Vert^2}{2} + \int_\Omega F_{\text{reg}}(u,v)dx + \dfrac{\sigma}{2}\|v-\overline{v}\|^2_*.
	\end{equation}
	Then we have
	\begin{lemma} \label{lem:omegastationary}
		Let $(u_\infty, v_\infty) \in \omega(u_0,v_0)$. Then $(u_\infty, v_\infty)$ is a critical point of $\eneg$.
	\end{lemma}
	\begin{proof}
Recalling the definition of the $\|\cdot\|_*$ norm, we calculate the first Fr\'{e}chet derivative of $\eneg$, namely,
	\[ \left\langle \eneg'(u,v), (h,k) \right\rangle = \varepsilon_u^2 \int_\Omega \nabla u \cdot \nabla h \: \mathrm{d}x + \varepsilon_v^2 \int_\Omega \nabla v \cdot \nabla k \: \mathrm{d}x + \int_\Omega \nabla F_\text{reg}(u,v) \cdot (h,k) \: \mathrm{d}x + \sigma \int_\Omega \mathcal{N}(v)k \: \mathrm{d}x,\]
	for $(u,v),(h,k) \in V_0^2$.
		Integrating by parts the first two terms and expanding the third one, we get
		\begin{multline*}
			\left\langle \eneg'(u_\infty,v_\infty), (h,k) \right\rangle =  \int_\Omega \left( -\varepsilon_u^2 \Delta u_\infty + \pd{F_\text{reg}}{u}(u_\infty, v_\infty) \right)h\: \mathrm{d}x \\ + \int_\Omega \left( -\varepsilon_v^2 \Delta v_\infty + \pd{F_\text{reg}}{v}(u_\infty, v_\infty) + \sigma \mathcal{N}(v_\infty) \right)k  \: \mathrm{d}x = 0,
		\end{multline*}
		since every point in the $\omega$-limit set is a stationary solution and lies in $Q_\epsilon$ for all $x \in \overline{\Omega}$ (we recall that $\overline{h} = \overline{k} = 0$ for all $(h,k) \in V_0^2$).
	\end{proof} \noindent
	We now show that $\eneg$ is twice continuously Fr\'{e}chet differentiable.
	\begin{lemma}
		The second Fr\'{e}chet derivative of $\eneg$ is well defined and $\eneg''(u,v) : V_0^2 \to \mathcal{L}(V_0^2, (V_0^2)^*)$ is given by
		\begin{multline*}
			\langle \eneg''(u,v)(w,z), (h,k) \rangle = \varepsilon_u^2 \int_\Omega \nabla w \cdot \nabla h \: \mathrm{d}x + \varepsilon_v^2 \int_\Omega \nabla z \cdot \nabla k \: \mathrm{d}x \\ + \int_\Omega (w,z) D^2 F_\text{reg}(u,v) \cdot (h,k) \: \mathrm{d}x + \sigma \int_\Omega \mathcal{N}(z)k \: \mathrm{d}x,
		\end{multline*}
		for all $(u,v), (w,z), (h,k) \in V_0^2$. Moreover, $\eneg''$ is continuous.
	\end{lemma}
	\begin{proof}
		We focus on the nonlocal term, since the computation of both the linear and nonlinear terms is straightforward (it suffices to expand up to order 1 the partial derivatives of $F_\text{reg}$). Moreover, since all derivatives of the regularized potential are uniformly bounded up to order three, then also continuity is an immediate consequence. \color{black} Let $H: V_0 \to V_0^*$ be defined as
		\[ \langle H(v), k \rangle = \sigma \ii{\Omega}{}{\mathcal{N}(v)k}{x}, \quad v,  k \in V_0.\]
		Notice that $H$ is linear and continuous, since the operator $\mathcal{N}$ is. Therefore, we have that arbitrary Fr\'{e}chet derivatives of $H$ exist, and they are all constant (and therefore continuous), namely:
		\[
		H': V_0 \to \mathcal{L}(V_0, V_0^*), \qquad H'(v) = H, \quad \forall\,v \in V_0.
		\]
		Moreover, notice that any variation of $H'$ vanishes, therefore $H'' \equiv 0$ (as an operator $H'': V_0 \to \mathcal{L}(V_0;\mathcal{L}(V_0,V_0^*))$), and all higher-order derivatives are zero as well. In fact, $H \in C^\infty(V_0,V_0^*)$. The proof is complete. \color{black}
	\end{proof} \noindent

Let us now introduce the linear operator $\mathcal{A}: V^2_0 \to (V^2_0)^*$ defined as follows (see Section \ref{sec:notation})
	\[
		\langle \mathcal{A}(u,v), (h,k) \rangle = \varepsilon_u^2 \langle Au, h\rangle + \varepsilon_v^2 \langle Av, k\rangle, \quad (u,v), (h,k) \in V_0^2.
	\]
	Notice that the bilinear form $a: V_0^2 \times V_0^2 \to \mathbb{R}$ defined by $a((u,v),(h,k)) := \langle \mathcal{A}(u,v), (h,k) \rangle$ is continuous and symmetric by the properties of $A$ and its inverse $\mathcal{N}$. Moreover, endowing $V_0$ with the classical $H^1$-seminorm,
	\[
	\langle \mathcal{A}(u,v), (u,v) \rangle = \varepsilon_u^2 \langle Au, u\rangle + \varepsilon_v^2 \langle Av, v\rangle = \varepsilon_u^2\|\nabla u \|^2 + \varepsilon_v^2\|\nabla v\|^2 \geq \min(\varepsilon_u^2,\varepsilon_v^2)\|(u,v)\|_{V_0^2}^2,
	\]
	for all $(u,v) \in V_0^2$, and thus coercivity holds as well. The corresponding operator $\mathcal{A}$ has therefore a nonempty resolvent set.	

The following lemmas are helpful in order to avoid \color{black} notational \color{black} ambiguities in the following.
	\begin{lemma} \label{lem:isomorphic}
		Let $X$ be a real Banach space. The space $(X^*)^2 = X^* \times X^*$ is isomorphic to $(X^2)^* = (X \times X)^*$. Moreover, the function $I : (X^*)^2 \to (X^2)^*$ acting as
		\[
			\langle I(L, M), (w,s) \rangle := \langle L, w \rangle + \langle M, s \rangle, \quad (L,M) \in (X^*)^2,\; (w,s) \in X^2
		\]
		is an isomorphism. Consequently, if we identify $H_0$ and its dual, the identification $H_0^2 \cong (H_0^2)^* \cong (H_0^*)^2$ is admissible.
	\end{lemma}
	\begin{proof}
		The linearity of $I$ is easily checked and follows from the linearity of its arguments. Assume that
		\[
			\langle I(L, M), (w,s) \rangle := \langle L, w \rangle + \langle M, s \rangle = 0
		\]
		for all $(w,s) \in X^2$. Testing on couples of functions in $X^2$ of the kind $(0,s)$, $(w,0)$, for arbitrary $w,s \in X$ immediately yields that both $L$ and $M$ must be the null operator in $X^*$. Therefore, $I$ is injective. Let now $P \in (X^2)^*$.
		Since, by linearity,
		\[ \langle P, (w,s) \rangle = \langle P, (w,0) \rangle + \langle P, (0,s) \rangle =: \langle P_w, w \rangle + \langle P_s, s \rangle, \quad (w,s) \in X^2, \]
		where the definitions of $P_w, P_s$ are clear from the equality above, if $P_w, P_s \in X^*$, then the operator $I$ is also surjective since the right hand side would equal $\langle I(P_w, P_s), (w,s) \rangle$. The linearity of $P_w, P_s$ comes from the linearity of $P$. As for continuity, recalling that we endow the product spaces with the respective standard Euclidean norms,
		\[
		\| P_w \|_{X^*} = \sup_{\|w\|_{X} = 1} |\langle P, (w,0) \rangle | \leq \sup_{\|w\|_{X} = 1} \|P\|_{(X^2)^*}\|(w,0)\|_{X^2} = \sup_{\|w\|_{X} = 1} \|P\|_{(X^2)^*}\|w\|_{X} = \|P\|_{(X^2)^*},
		\]
		proving that $P_w \in X^*$. An analogous proof works for $P_s$. Thus, the operator $I$ is invertible and we are only left to prove its continuity. Indeed,
		\[
		\|I(L,M)\|_{(X^2)^*} = \sup_{\|(w,s)\|_{X^2} = 1} |\langle I(L,M), (w,s) \rangle| = \sup_{\|(w,s)\|_{X^2} = 1} | \langle L, w \rangle | + |\langle M, s \rangle|
		\]
		by the triangle inequality, and furthermore, by continuity of $L$ and $M$
		\[
		\begin{split}
			\sup_{\|(w,s)\|_{X^2} = 1} | \langle L, w \rangle | + |\langle M, s \rangle| & \leq  \sup_{\|(w,s)\|_{X^2} = 1} \|L\|_{X^*}\|w\|_{X} + \|M\|_{X^*}\|s\|_{X}\\ &  \leq \max \left(\|L\|_{X^*},\|M\|_{X^*}\right)  \sup_{\|(w,s)\|_{X^2} = 1} \left(\|w\|_{X} + \|s\|_{X}\right),
		\end{split}
		\]
		and, since all norms are equivalent in $\mathbb{R}^2$, we deduce that there exists $C > 0$ independent of $L,M,w$ and $s$ such that
		\[	\|I(L,M)\|_{(X^2)^*}  \leq  \max \left(\|L\|_{X^*},\|M\|_{X^*}\right)  \sup_{\|(w,s)\|_{X^2} = 1} \left(\|w\|_{X} + \|s\|_{X}\right) \leq C(\|L\|_{X^*}^2 + \|M\|_{X^*}^2)^\frac{1}{2},\]
		and the first part of the statement is proved. As for the second part of the statement, we consider the case $X = H_0$. First of all notice that the Hilbert triplet setting lets us identify the set $H_0$ with its dual space $H^*_0$. Observe that the Riesz isometry $R_2$ between $H_0^2$ and its dual is defined by
		\[ \langle R_2(w,s), (h,k) \rangle := \langle Rw, h \rangle + \langle Rs, k \rangle, \quad (w,s), \: (h,k) \in H_0^2,\]
		where $R$ denotes the Riesz map between $H_0$ and its dual, then we can also identify $H_0^2$ and its dual, and moreover, in this case, after identification,
		\[
			R_2(w,s) = I(Rw, Rs) \to (w,s) = I(w,s),
		\]
		and thus we also deduce that $I$ is the identity map (as expected).
		\end{proof}
	\begin{lemma} \label{lem:domain}
		Let $\mathcal{A}_2$ be the part of $\mathcal{A}$ in $H_0^2$, namely
			\[ \mathcal{A}_2 : \mathcal{D}(\mathcal{A}_2) \to H_0^2,\]
		where the domain $\mathcal{D}(\mathcal{A}_2) := \{ (w,s) \in V_0^2: \mathcal{A}(w,s) \in H_0^2\}.$ Furthermore, let \[ A_2 : \mathcal{E}(A_2) \to H_0\] denote the part in $H_0$ of the operator $A$ defined in Section \ref{sec:notation} with domain $\mathcal{E}(A_2)$. Then $\mathcal{D}(\mathcal{A}_2) = \mathcal{E}(A_2) \times \mathcal{E}(A_2)$.
	\end{lemma}
	\begin{proof}
		The inclusion $\mathcal{E}(A_2) \times \mathcal{E}(A_2) \subset \mathcal{D}(\mathcal{A}_2)$ holds, since if $(w,s) \in \mathcal{E}(A_2) \times \mathcal{E}(A_2)$, then by definition $w,s \in V_0$ are such that $Aw, As \in H_0$. Therefore, owing to Lemma \ref{lem:isomorphic}, the equality
		\[ \mathcal{A}(w,s) = I(Aw, As) = (Aw, As) \in H_0^2 \]
		holds in the $H_0^2$-sense. Conversely, let now $(w,s) \in \mathcal{D}(\mathcal{A}_2)$. Then,
		\[ \mathcal{A}(w,s) = I(Aw, As) \in H_0^2 \Rightarrow Aw, As \in H_0 \Rightarrow Aw, As \in \mathcal{E}(A_2),\]
		and the statement is proved.
	\end{proof}
	\begin{remark} \label{rem:domain}
		As customary, from Lemma \ref{lem:domain} we also know that
		\[ \mathcal{D}(\mathcal{A}_2) = \left\{ (u,v) \in H^2(\Omega)^2: \pd{u}{\mathbf{n}} = \pd{v}{\mathbf{n}} = 0 \text{ a.e. on } \partial \Omega \right\}.\]
	\end{remark}
	\begin{lemma} \label{lem:kernel}
		Let $(u,v) \in V_0^2$ and let $(w,z) \in \ker \eneg''(u,v)$. Then $(w,z) \in \mathcal{D}(\mathcal{A}_2)$.
	\end{lemma}
	\begin{proof}
		If $(w,z)$ belongs to the kernel of $\eneg''$, then
		\begin{multline*}
			\langle \eneg''(u,v)(w,z), (h,k) \rangle = \varepsilon_u^2 \int_\Omega \nabla w \cdot \nabla h \: \mathrm{d}x + \varepsilon_v^2 \int_\Omega \nabla z \cdot \nabla k \: \mathrm{d}x \\ + \int_\Omega (w,z) D^2 F_\text{reg}(u,v) \cdot (h,k) \: \mathrm{d}x + \sigma \int_\Omega \mathcal{N}(z)k \: \mathrm{d}x = 0,
		\end{multline*}
		for every choice of $(h,k) \in V_0^2$. Expanding the third integral and separating variables we obtain the following equality (we omit the dependence of the partial derivatives on $(u,v)$ for the sake of clarity):
		\[
			\left \langle \varepsilon_u^2 Aw + \left( \dfrac{\partial^2 F_\text{reg}}{\partial u^2}+ \dfrac{\partial^2 F_\text{reg}}{\partial u \partial v} \right)w, h \right\rangle = \left \langle - \varepsilon_v^2 Az -\left( \dfrac{\partial^2 F_\text{reg}}{\partial v^2}+ \dfrac{\partial^2 F_\text{reg}}{\partial u \partial v} \right)z - \sigma \mathcal{N}z, k \right \rangle,
		\]
		for every $(h,k) \in V_0^2$. Testing on  $(h,0)$, $(0,k) \in V_0^2$ we infer that both functionals must equal the null operator in $V_0^*$, thus there hold
		\[
			\begin{cases}
				Aw = -\dfrac{1}{\varepsilon_u^2}\left( \dfrac{\partial^2 F_\text{reg}}{\partial u^2}+ \dfrac{\partial^2 F_\text{reg}}{\partial u \partial v} \right)w, \\[0.3cm]
				Az = -\dfrac{1}{\varepsilon_v^2}\left[\left( \dfrac{\partial^2 F_\text{reg}}{\partial v^2}+ \dfrac{\partial^2 F_\text{reg}}{\partial u \partial v} \right)z - \sigma \mathcal{N}z\right],
			\end{cases}
		\]	
		where the equalities are to be intended in the $V^*_0$-sense. However, by uniform boundedness of the derivatives, and since $\mathcal{N}z \in V_0 \subset H_0$, it is easy to notice that both right hand sides are well defined in $H_0$ (we recall that $(w,z) \in V_0^2$ as well). The statement is proved.
	\end{proof} \noindent
	In more than one spatial dimension the set $\mathcal{E}$ is usually a continuum. Therefore, one cannot simply take advantage of the fact that we are dealing with a gradient system. An appropriate tool in this case is a method based on the so-called {\L}ojasiewicz-Simon inequality which is given, in our case, by
	\begin{proposition} \label{prop:lsineq}
		Let $(u_0,v_0) \in Z$.
Suppose that $(u_\infty, v_\infty) \in \omega(u_0,v_0)$ and \color{black} assume that the singular part $\hat{S}$ is real analytic (see Remark \ref{genpot})\color{black}. Then, there exist constants $\theta \in (0,\frac{1}{2})$ and $C, \varpi > 0$ such that
		\[ |\eneg(w,z) - \eneg(u_\infty, v_\infty)|^{1-\theta} \leq C\|\eneg'(w,z)\|_{(V_0^2)^*} ,\]
		provided that $\|(w,z) - (u_\infty, v_\infty)\|_{V^2} \leq \varpi$.
	\end{proposition}
\begin{remark}
We remark that a necessary requirement is the real analyticity of the nonlinearities. Indeed, we recall that even for $\mathcal{C}^\infty$ nonlinearities there are counterexamples (see, for instance, \cite{PS}). Also, we recall that it might happen that $\theta=\frac{1}{2}$ but this depends on the possible hyperbolic nature of the stationary states (see, e.g., \cite[Cor. 3.13]{Chill03}).
\end{remark}
	\begin{proof}
		By Lemma \ref{lem:domain} and Remark \ref{rem:domain}, as well as Sobolev embedding, we have that $\mathcal{D}(\mathcal{A}_2) \hookrightarrow L^\infty(\Omega)^2$. Note that both operators $\mathcal{A}$ and $\mathcal{A}_2$ have compact resolvents in $(V_0^2)^*$ and $H_0^2$, respectively. Moreover, it follows by Proposition \ref{prop:maxprinc} that any $(u_\infty, v_\infty) \in \omega(u_0,v_0)$ is uniformly bounded, and thus the operators $\eneg''$ and its restriction to $\mathcal{D}(\mathcal{A}_2)$ are bounded perturbations of $\mathcal{A}$ and $\mathcal{A}_2$, respectively. Thus, arguing as in \cite[Proposition 6.6]{Chill06}, the space $\ker \eneg''(u_\infty,v_\infty)$ is a finite-dimensional subspace of $\mathcal{D}(\mathcal{A}_2)$ (owing also to Lemma \ref{lem:kernel}). Also, the range $\rg \eneg''(u_\infty,v_\infty)$ is closed in $(V_0^2)^*$ and the range $\rg \restr{\eneg''(u_\infty,v_\infty)}{\mathcal{D}(\mathcal{A}_2)}$ is closed in $H_0^2$. We want to apply \cite[Corollary 3.11]{Chill03}. First of all, notice that, by Lemma \ref{lem:omegastationary}, any point in the $\omega$-limit set is a stationary point of the energy functional. Adopting the same notation as in \cite{Chill03}, we set $X := \mathcal{D}(\mathcal{A}_2) \subset V_0^2$, $Y := H^2_0$, $W := (V_0^2)^*$. Moreover, let us denote by $\Pi^* : (V_0^2)^* \to (V_0^2)^*$ the orthogonal projector onto $\ker \eneg(u_\infty,v_\infty)$. By Lemma \ref{lem:kernel}, $\Pi^*(X) \subset X$, and, furthermore, $\Pi^*(Y) \subset Y$. Moreover, the energy functional has the form
		\[
		\eneg(u,v) = \dfrac{1}{2} \langle \mathcal{A}(u,v), (u,v) \rangle + \dfrac{\sigma}{2} \langle v, \mathcal{N}v \rangle + \int_\Omega F_\text{reg}(u,v) \: \mathrm{d}x,
		\]
		and since $F$ is analytic in $Q_\epsilon$, then the Fr\'{e}chet derivative $\eneg'$ is real analytic (as an operator between Banach spaces) in a neighborhood of $(u_\infty, v_\infty)$ in $X$. We conclude the result using \cite[Corollary 3.11]{Chill03}.
	\end{proof} \noindent
	We can now state and prove the convergence to a single stationary state.
	\begin{theorem}
		Let $\hat{S}$ be real analytic in $(-1,1)$. Let $(u_0,v_0) \in V^2_0$ such that $F(u_0,v_0) \in L^1(\Omega)$. Consider the trajectory $(u(t),v(t))=\mathcal{S}(t)(u_0,v_0)$. Then there exists $(u_\infty, v_\infty) \in \mathcal{E}$ such that
		\[
		\lim_{t \to +\infty} (u(t), v(t)) = (u_\infty, v_\infty), \qquad \text{ in } H^{2r}(\Omega)^2,
		\]
		for all $r \in (0,1)$. Moreover, there exists $C>0$, depending also on $\theta$, such that, for all $t\geq 0$,
\begin{equation}
\label{eq:rateconv}
\Vert (u(t),v(t)) - (u_\infty,v_\infty)\Vert_{(V^*_0)^2} \leq C(1+t)^{-\frac{\theta}{1-2\theta}}.
\end{equation}
\end{theorem}
	\begin{proof}
		The LaSalle's invariance principle and Proposition \ref{prop:maxprinc} entail that
		\[
			\restr{\eneg(u,v)}{\omega(u_0,v_0)} = \restr{\ene(u,v)}{\omega(u_0,v_0)} \equiv \widetilde{\Psi}_\infty,
		\]
		for some $\widetilde{\Psi}_\infty \in \mathbb{R}$. By compactness of $\omega(u_0,v_0)$, we can consider an open cover $U_2$ formed by a finite number of sufficiently small balls, namely
		\[
		\omega(u_0,v_0) \subset U_2 := \bigcup_{i = 1}^N  B_{\varpi_i}((u_i, v_i)),
		\]
		for some $N \in \mathbb{N}$, $(u_i, v_i) \in V_0^2$ for $i = 1, \dots, N$ and radii $\varpi_i > 0$ satisfying
		\[ \max_{1 \leq i \leq N} \varpi_i \leq \varpi,\]
		where $\varpi$ is the one appearing in Proposition \ref{prop:lsineq}. Notice that Proposition \ref{prop:lsineq} holds in each of the balls, and since they are a finite number, we can extract uniform constants $C > 0, \theta \in (0,\frac{1}{2}]$ such that
		\begin{equation}
\label{eq:LSin}
			|\eneg(w,z) - \widetilde{\Psi}_\infty|^{1-\theta} \leq C\|\eneg'(w,z)\|_{(V_0^2)^*} , \quad (w,z) \in U.
		\end{equation}
		Observe that also $U_2$ attracts the trajectory of the dynamical system. Thus $(u(t),v(t)) \in U_2$ for every $t \geq t_2$. Analogously, if $t_1$ is such that $(u(t),v(t)) \in U_1$ for every $t \geq t_1$, we consider the trajectory starting from the time instant $t^\sharp := \max(t_1, t_2)$, so that $(u(t), v(t)) \in U_1 \cap U_2 \supset \omega(u_0,v_0)$ for all $t \geq t^\sharp$. In this way,
		\[ (u(t,x), v(t,x)) \in Q_\epsilon, \qquad \forall t \geq t^\sharp, \; \forall\, x \in \Omega.\]
		Let us consider the functional $H : [t^\sharp, +\infty) \to \mathbb{R}$ defined by
		\[
			H(t) = (\eneg(u(t),v(t)) - \widetilde{\Psi}_\infty)^\theta
		\]
		Then, a straightforward computation, jointly with the differential form of Lemma \ref{lem:identity} and \ref{eq:LSin}, entail that
		\begin{equation} \label{eq:conv1}
					-\td{}{t}H(t) = -\theta (\eneg(u(t),v(t)) - \widetilde{\Psi}_\infty)^{\theta-1} \td{}{t}\eneg(u(t),v(t)) \geq C\theta \dfrac{\|\nabla \mu(t) \|^2 + \|\nabla \widetilde{\varphi} (t) \|^2}{\|\eneg'(u(t),v(t))\|_{(V_0^2)^*}}.
		\end{equation}
		Next, we consider the quantity $\|\eneg'(u(t),v(t))\|_{(V_0^2)^*}$. By definition,
		\[
		\begin{split}
					\|\eneg'(u(t),v(t))\|_{(V_0^2)^*}  = \sup_{\|(h,k)\|_{V_0^2} = 1} |\langle \eneg'(u(t),v(t)), (h,k) \rangle|.
		\end{split}
		\]
		Recalling the computations in Lemma \ref{lem:omegastationary} and the fact that $(u(t),v(t)) \in U_1$, the right hand side equals \color{black}
		\[
		\begin{split}
			\sup_{\|(h,k)\|_{V_0^2} = 1} & \left| \int_\Omega \left( -\varepsilon_u^2 \Delta u(t) + \pd{F}{u}(u(t),v(t)) \right) h\: \mathrm{d}x + \int_\Omega \left( -\varepsilon_v^2 \Delta v(t) + \pd{F}{v}(u(t),v(t)) + \sigma \mathcal{N}(v(t)) \right)k  \: \mathrm{d}x\right| \\
			& = \sup_{\|(h,k)\|_{V_0^2} = 1} \left| \int_\Omega \mu(t) h \: \mathrm{d}x + \int_\Omega \widetilde{\varphi}(t) k \: \mathrm{d}x \right| \\
			& = \sup_{\|(h,k)\|_{V_0^2} = 1} \left| \int_\Omega \left(\mu(t) - \overline{\mu(t)}\right) h \: \mathrm{d}x + \int_\Omega \left(\widetilde{\varphi}(t) - \overline{\widetilde{\varphi}(t)}\right) k \: \mathrm{d}x \right| \\
			& \leq C\sup_{\|(h,k)\|_{V_0^2} = 1} \| \nabla \mu(t) \| \|h\| + \| \nabla \widetilde{\varphi}(t)\| \| k \| \\
			& \leq C\max (\| \nabla \mu(t) \|, \| \nabla \widetilde{\varphi}(t)\|) \sup_{\|(h,k)\|_{V_0^2} = 1} (\|h\| + \|k\|) \\
			& \leq C \left( \| \nabla \mu(t) \| + \| \nabla \widetilde{\varphi}(t) \| \right),
		\end{split}
		\] \color{black}
		for some $C > 0$, possibly different from the one above. We also used the Poincar\'{e}-Wirtinger and the triangle inequalities, as well as the equivalence of all norms in $\mathbb{R}^2$. Therefore from \eqref{eq:conv1} we deduce
		\[
		 -\td{}{t}H(t) \geq C \dfrac{\|\nabla \mu(t) \|^2 + \|\nabla \widetilde{\varphi} (t) \|^2}{\| \nabla \mu(t) \|+ \| \nabla \widetilde{\varphi}(t) \| } \geq \dfrac{C}{2} \left( \| \nabla \mu(t) \| + \| \nabla \widetilde{\varphi}(t) \| \right),
		\]
		and an integration on $[t^\sharp, +\infty)$ yields that $\nabla \mu$, $\nabla \widetilde{\varphi} \in L^1([t^\sharp,+\infty), H_0)$. By comparison in the evolution equations, we also have $\partial_t u, \partial_t v \in L^1([t^\sharp,+\infty), V^*_0)$. Therefore, the limit
		\[ \lim_{t \to +\infty} (u(t),v(t)) = (u_\infty, v_\infty) \]
		holds in the $V_0^*$-sense and, by compactness of the $\omega$-limit set, also in the $H^{2r}$-sense, for all $r \in (0,1)$. The fact that $(u_\infty, v_\infty)\in \mathcal{E}$ follows from the characterization of the stationary points of $\mathcal{S}$. \\
We have $\widetilde{\Psi}_\infty = \eneg(u_\infty,v_\infty)$. Observe that, thanks to \eqref{eq:coupledCHequi}, we have
\begin{align}
\Vert (u(t),v(t)) - (u_\infty,v_\infty)\Vert_{(V^*_0)^2} &\leq C \int_t^{+\infty} \Vert (u_t(\tau),v_t(\tau))\Vert_{(V^*_0)^2} \:\mathrm{d}\tau \nonumber\\
&\leq C \int_t^{+\infty} \Vert (\nabla \mu(\tau), \nabla\widetilde{\varphi}(\tau))\Vert\:\mathrm{d}\tau  \leq C H(t), \label{eq:ratein}
\end{align}
for all $t\geq t^\sharp$. From \eqref{eq:eneriden} and \eqref{eq:conv1} we deduce
\begin{equation*}
\td{}{t}H(t) + \theta H(t)^{\frac{\theta-1}{\theta}} \td{}{t}\eneg(u(t),v(t)) = \td{}{t}H(t) + \theta H(t)^{\frac{\theta-1}{\theta}}\left(\| \nabla \mu(t) \|^2 + \| \nabla \widetilde{\varphi}(t) \|^2\right)  =0.
\end{equation*}
On the other hand, using again \eqref{eq:LSin}, we get
$$
C H(t)^{{\frac{2(1-\theta)}{\theta}}} \leq \left(\| \nabla \mu(t) \|^2 + \| \nabla \widetilde{\varphi}(t) \|^2\right) , \qquad \forall\,t\geq t^\sharp,
$$
so that
\begin{equation}
\label{eq:conv2}
\td{}{t}H(t) + \theta C H(t)^{\frac{1-\theta}{\theta}}  \leq 0, \qquad \forall\,t\geq t^\sharp.
\end{equation}
From \eqref{eq:ratein} and \eqref{eq:conv2} we deduce the wanted estimate rate \eqref{eq:rateconv} (see \cite[Cor. 6.3.3]{SZ}).
\end{proof}

	\begin{remark} \label{rem:nonzeromean}
		If the initial conditions do not have null spatial average, that is, for instance, $\overline{u_0} = m_1$ and $\overline{v_0} = m_2$ for some $m_1, m_2 \in (-1,1)$, then we can always reformulate the problem in order to have zero mean and argue as above. Indeed, setting $\tilde{u} := u - m_1$ and $\tilde{v} := v-m_2$, we have $\tilde{u}, \tilde{v} \in V_0$, and the pair $(\tilde{u},\tilde{v})$ solves the \color{black}problem\color{black}
		\begin{equation*}
			\begin{cases}
				\pd{\tilde{u}}{t} = \Delta \mu & \quad \text{in } \Omega \times (0,T),\\
				\mu = - \varepsilon_u^2 \Delta \tilde{u} + \pd{F}{u}(\tilde{u} + m_1, \tilde{v}+m_2), & \quad \text{in } \Omega \times (0,T),\\
				\pd{\tilde{v}}{t} = \Delta \widetilde{\varphi} & \quad \text{in } \Omega \times (0,T),\\
				\widetilde{\varphi} = - \varepsilon_v^2 \Delta \tilde{v} + \pd{F}{v}(\tilde{u}+m_1,
				\tilde{v}+m_2) + \sigma \mathcal{N}\tilde{v} & \quad \text{in }  \Omega \times (0,T),\\
				\pd{\tilde{u}}{\mathbf{n}} = \pd{\tilde{v}}{\mathbf{n}} = 0 & \quad \text{on } \partial \Omega \times (0,T),\\[0.2cm]
				\pd{\mu}{\mathbf{n}} = \pd{\varphi}{\mathbf{n}} = 0 & \quad \text{on } \partial \Omega \times (0,T),\\[0.2cm]
				u(\cdot, 0) = u_0 & \quad \text{in } \Omega,\\
				v(\cdot, 0) = v_0 & \quad \text{in } \Omega.
			\end{cases}
		\end{equation*}
		Therefore, replacing the nonlinearity $F$ with $\widetilde{F}(s_1, s_2) := F(s_1+m_1, s_2 +m_2)$, we recover the structure of Problem \eqref{eq:coupledCHequi}. Notice that this change does not affect the regularity of the nonlinear term.
	\end{remark}

\begin{remark}
On account of Proposition \ref{prop:reg3}, using the general theory, one can prove that the dynamical system $(Z,\mathcal{S}(t))$ has a connected global attractor bounded in $H^4(\Omega)^2$ if $d=2$ (see Remark \ref{2REG}) or $W^{2,6}(\Omega)^2$ if $d=3$ (see \cite[Thm. 6.1]{GGM}). This attractor coincides with the unstable manifold of $\mathcal{E}$
since $\mathcal{S}(t)$ has a Lyapunov function on $Z$. If $d=2$ then the separation property allows us to establish the existence of an exponential attractor  (see \cite[Thm. 6.1]{GGM} and related remarks).
\end{remark}

\section{Concluding remarks and future issues}
\label{sec:conclusion}

The analysis of regularity and its implications in the off-critical case remains an open issue. As we pointed out in Section \ref{sec:regularity}, because $\overline{\dq{v}}$  is no longer zero (see \eqref{offcriticalquot}) we are unable to carry out the proof of the crucial Proposition \ref{prop:reg1} as we did. A different strategy might be required. Taking the hydrodynamic effects into account could be another challenging issue (see, for instance, \cite{Miranville2} and its references) as well as replacing the standard Cahn-Hilliard equation with its nonlocal counterpart (see \cite{GGG} and references therein). We also recall that some models of surfactants are represented by coupled Cahn-Hilliard equations (see \cite{KK}, cf. also \cite{Yang18} and its references for the numerical approximations) possibly with hydrodynamic effects (see, e.g., \cite{EDAT} and references therein).
The present approach could be extended to these models. Indeed they are characterized by regular potentials so that one cannot ensure that the local concentrations take their values in the physical range. As a consequence, one cannot guarantee that the total free energy is bounded from below. Instead, we believe that taking mixing entropies as we did here can led us to establish physically
meaningful theoretical results. This choice might also help to design alternative numerical schemes (compare with the penalization argument introduced in \cite{XCY}). \color{black} Concerning the longtime behavior, the existence of an exponential attractor, which entails the finite-dimensionality of the global attractor, cannot be extended easily to the case $d=3$ since no (global) separation property is known. However, one might try to use the argument devised in \cite[Sec.5]{Miranville3}, a sort of ``local'' separation property.
\color{black}

\bigskip
\noindent
{\bf Acknowledgment}. The authors thank the reviewers for their useful comments. The second author is a member of Gruppo Nazionale per l'Analisi Matematica, la Probabilit\`{a} e le loro Applicazioni (GNAMPA), Istituto Nazionale di Alta Matematica (INdAM).

	\printbibliography[title=Bibliography]
\end{document}